\documentclass[leqno,12pt]{article}
\usepackage{color}
\usepackage{ifthen}
\usepackage{soul}
\usepackage{latexsym}
\usepackage{amsthm,amssymb}
\usepackage{amsbsy,amsfonts,amsmath}

\usepackage{a4}
\usepackage{amsthm,amssymb}
\usepackage{amsbsy,amsfonts,amsmath}

\usepackage{latexsym}

\newcommand{\al}{{\alpha}}

\newcommand{\cR }{\mathcal{R}}

\newcommand{\dg }{{\rm{deg}}}

\newcommand{\End }{\mathrm{End}}
\newcommand{\id }{\mathrm{id}}

\newcommand{\ndC }{\mathbb{C}}
\newcommand{\ndN }{\mathbb{N}}
\newcommand{\ndQ }{\mathbb{Q}}
\newcommand{\ndR }{\mathbb{R}}
\newcommand{\ndZ }{\mathbb{Z}}

\newcommand{\ndA }{\mathbb{A}}

\newcommand{\psl }{\mathfrak{psl}}

\newcommand{\srs}{\mathcal{R}}
\newcommand{\rs}{R}
\newcommand{\s}{\sigma }
\newcommand{\GL }{\mathrm{GL}}

\newcommand{\is }{A}

\newcommand{\sa}{\triangleright}
\newcommand{\vecm}{\mathbf{m}}
\newcommand{\Wmon}{\widetilde{W}}
\newcommand{\sig }[2]{#1\sa #2}

\newcommand{\te}{\tilde{e}}
\newcommand{\ts}{\tilde{s}}

\newcommand{\tr}{\triangleright}

\newtheorem{theorem}{Theorem}[section]
\newtheorem{proposition}[theorem]{Proposition}
\newtheorem{corollary}[theorem]{Corollary}
\newtheorem{lemma}[theorem]{Lemma}

\theoremstyle{definition}
\newtheorem{defin}[theorem]{Definition}
\newtheorem{notation}[theorem]{Notation}
\newtheorem{example}[theorem]{Example}
\newtheorem{remark}[theorem]{Remark}

\title{Iwahori-Hecke type algebras \\
associated with the Lie superalgebras \\
$A(m,n)$, $B(m,n)$, $C(n)$ and $D(m,n)$}

\author{Hiroyuki Yamane}
\date{}

\begin{document}

\maketitle

\begin{abstract}
In this paper we give Iwahori-Hecke type algebras
$H_q(\mathfrak{g})$
associated with the Lie superalgebras 
$\mathfrak{g}=A(m,n)$, $B(m,n)$, $C(n)$ and $D(m,n)$.
We classify the irreducible representations
of $H_q(\mathfrak{g})$
for generic $q$.
\end{abstract}

\section*{Introduction}
Recently, motivated by a question posed by V.~Serganova~\cite{a-Serg96}
and study of the Weyl groupoids \cite{a-Heck06b}\cite{a-Heck04c} associated with
Nichols algebras \cite{a-AndrSchn98}\cite{a-AndrSchn05p}
including generalizations
of quantum groups, 
I.~Heckenberger and the author \cite{HY} introduced a notion of 
`Coxeter groupoids' (in fact they can be defined as semigroups), and showed that a Matsumoto-type theorem 
holds for the groupoids, so they have the solvable word problem.
We mention that the Coxeter groupoid associated with
the affine Lie superalgebra $D^{(1)}(2,1;x)$
was used in the study \cite{HSTY},
where Drinfeld second realizations of $U_q(D^{(1)}(2,1;x))$ was analized 
by
physical motivation in recent study of AdS/CFT correspondence.

It would be able to be said that one of the main purposes at present of
the representation theory is to study the Kazhdan-Lusztig polynomials
(cf. \cite[7.9]{Hu})
and their versions.
The polynomials are defined by using
the standard and canonical bases of the Iwahori-Hecke algebras. 
The existence of those bases is closely related to
the Matsumoto theorem of the Coxeter groups. So it would be natural
to ask what to be the Iwahori-Hecke algebras of the 
Coxeter groupoids. In this paper, we give a tentative answer to
this question for the Coxeter groupoids $W$ associated with
the Lie superalgebras $\mathfrak{g}=A(m,n)$, $B(m,n)$, $C(n)$ and $D(m,n)$. 
We introduce the Iwahori-Hecke type
algebra $H_q(\mathfrak{g})$ (in the text, it is also denoted by $H_q(W)$) 
as $q$-analogue of the 
semigroup algebra $\mathbb{C}W/\mathbb{C}0$, where $0$ is the zero element of $W$. 
We also show that 
if $q$ is nonzero and not any root of unity, 
$H_q(\mathfrak{g})$ is semisimple and 
there exists a natural one-to-one correspondence between 
the equivalence classes of the irreducible representations 
of $H_q(\mathfrak{g})$ 
and those of the Iwahori-Hecke algebra
$H_q(W_0)$
associated with the Weyl group
$W_0$ of
the Lie algebra $\mathfrak{g}(0)$
obtained as  
the even part of $\mathfrak{g}=\mathfrak{g}(0)\oplus \mathfrak{g}(1)$. 

Until now, 
no relation has been achieved between the groupoids
treated in \cite{SV} and this paper. 

This paper is composed of the two sections.
Main results and their proofs are given in Section~2.
Results of \cite{HY} used in Section~2 are introduced 
in Section~1. 

The author thanks to the referee for careful reading and valuable comments,
which encourage him so much to make future study. 

\section{Preliminary---Matsumoto-type theorem of Coxeter groupoids}
\label{section:preliminary}

This section is preliminary. Here we collect the results which have already 
been given in \cite{HY}
and will be used in the next section.

\subsection{Semigroups and Monoids}

Let $K$ be a non-empty set. Assume that $K$ has a product map
$K\times K\rightarrow K$, $(x,y)\mapsto xy$.
We call $K$ a {\it{semigroup}} if $(xy)z=x(yz)$ for $\forall$ $x,y,z\in K$.
We call $K$ a {\it{monoid}} if $K$ is a semigroup and
there exists a unit $1\in K$, that is, 
$1x=x1=x$ for all $x\in K$.
\subsection{Free semigroup  $F_{-1}(N)$ and Free monoid  $F_0(N)$}
Let $N$ be a non-empty set. Let $F_{-1}(N)$ be the set of all the finite sequences of 
elements of $N$, that is 
$$
F_{-1}(N):=\coprod_{n=1}^\infty N^n=\{(h_1,\ldots,h_n)|n\in\ndN, h_i\in N\}.
$$ We regard $F_{-1}(N)$ as the semigroup by 
$$
(h_1,\ldots,h_m)(h_{m+1},\ldots, h_{m+n})
=(h_1,\ldots, h_m,h_{m+1},\ldots, h_{m+n}).
$$ Then we call $F_{-1}(N)$ a {\it{free semigroup}}.
Let $F_0(N)$ be the semigroup obtained by adding the unit $1$, that is,  
$F_0(N):=\{1\}\cup F_{-1}(N)$, $1\notin F_{-1}(N)$,
and $1x=x1=x$ for all $x\in F_0(N)$.
\subsection{Semigroup generated by the generators and 
and defined by the relations}
Let
$Q=\{(x_j,y_j)|j\in J\}$ be a subset of $F_{-1}(N)\times F_{-1}(N)$,
where $J$ is an index set.
For $g_1$, $g_2\in F_{-1}(N)$,
we write 
$g_1 \sim_1 g_2$
if there exist $j\in J$ and $(f_1, f_2)\in F_0(N)
\times F_0(N)$ such that either of the following (i), (ii), (iii) holds. 
\newline\par
(i) $g_1=f_1xf_2\ne g_2=f_1yf_2$.\par
(ii) $g_1=f_1yf_2\ne g_2=f_1xf_2$.\par
(iii) $g_1=g_2=f_1xf_2=f_1yf_2$.
\newline\newline
For $g$, $g^\prime\in F_{-1}(N)$,
we write 
$g \sim g^\prime$
if  
$g=g^\prime$ or
there exists $r\in \ndN$ and $g_1, \ldots , g_r \in F_{-1}(N)$
such that
$g_1=g$, $g_r=g^\prime$,
and $g_i \sim_1 g_{i+1}$ for $1\leq i\leq r-1$.
Then $F_{-1}(N)/\!\!\sim$
can be regarded as a semigroup by the product $[g][g^\prime]=[gg^\prime]$,
where for $g\in F_{-1}(N)$, we denote $[g]:=\{g^\prime|g^\prime\sim g\}\in
F_{-1}(N)/\!\!\sim$.
We call 
$F_{-1}(N)/\!\!\sim$ 
{\it{the semigroup generated by $N$ 
and defined by the relations $x_j=y_j$
($j\in J$). }}
When there is no fear of misunderstanding,
we also denote $[g]$ by its representative 
$g$ by abuse of notation.
\subsection{Free group  $F_1(N)$ and Involutive free group  $F_2(N)$} Let $N$ be a set. Let $N^{-1}$ be a copy of $N$
so that the bijective map $N\rightarrow N^{-1}$, $x\mapsto x^{-1}$, is given.
Let $F_1(N)$ be the semigroup generated by
$$
\{e\}\cup N \cup N^{-1} \quad \text{(disjoint union)}
$$ and defined by the relations
$$
ee=ex=xe=e,\quad xx^{-1}=x^{-1}x=e
\quad \text{for $\forall x\in N$.} 
$$ 
We call $F_1(N)$ the {\it{free group}} over $N$.

Let $F_2(N)$ be the semigroup generated by 
$$
\{e\}\cup N \quad \text{(disjoint union)}
$$ and defined by the relations
$$
ee=ex=xe=e,\quad x^2=e
\quad \text{for $\forall x\in N$.} 
$$ We call $F_2(N)$ the {\it{involutive free group}} over $N$.
Note that $F_2(N)$ can be identified with the quotient group of $F_1(N)$ in the natural sense:
$$
F_2(N)=F_1(N)/\{g_1x_1^2g_1^{-1}\cdots g_rx_r^2g_r^{-1}
|r\in\ndN, x_i\in N, g_i\in F_1(N)\}.
$$
\subsection{Action $\sa$ of $F_2(N)$ on $A$}
Let $N$ and $A$ be non-empty sets. An action $\sa$ of $F_2(N)$ on $A$ is
a map 
$$\sa :F_2(N)\times A \to A
$$ such that 
$$
e\sa a=a,\quad
g\sa (h\sa a)=(gh)\sa a \quad
\text{for 
$\forall g,\forall h\in F_2(N)$, $\forall a\in A$.}
$$ Note that
$n\sa (n\sa a)=a$ for all
$n\in N,a\in A$.

For $n$, $n^\prime\in N$ and $a\in A$, define
\begin{align*}
\Theta (n,n';a):=&\{(nn')^m\sa a,\,
(n'n)^m\sa a\,|\,m\in \ndN\cup\{0\}\}.
\end{align*}
Let $$
\theta (n,n';a):=|\Theta (n,n';a)|.$$
This is the cardinality of $\Theta (n,n';a)$,
which is either in $\ndN $ or is $\infty $.
One obviously has $\Theta (n,n';a)=\Theta (n',n;a)$
and $\Theta (n,n';n\sa a)=n\sa \Theta (n',n;a)$.
\newline\newline
Let $a_0:=a$, $b_0:=a$, and define recursively
$a_{m+1}:=n\sa b_m$, $b_{m+1}:=n'\sa a_m$ for all $m\in \ndN\cup\{0\}$. That is:
\begin{align*}
& b_0:=a, \quad a_1:=n \sa a, \quad b_2:=n' \sa n \sa a, \quad a_3:=n \sa n' \sa n \sa a,\ldots \\
& a_0:=a, \quad b_1:=n' \sa a, \quad a_2:=n \sa n' \sa a, \quad b_3:=n' \sa n \sa n' \sa a,\ldots
\end{align*} 
Then we have
\begin{align*}
 \theta (n,n';a)=
 \begin{cases}
  \infty & \text{if $a_m\not=b_m$ for all $m\in \ndN $,}\\
  \min \{m\in \ndN \,|\,a_m=b_m\} & \text{otherwise.}
 \end{cases}
\end{align*}

\subsection{Coxeter groupoids}

\begin{defin}\cite[Definition~1]{HY}\label{df:Cox}
Let $N$ and $A$ be non-empty sets.
Let
$\sa $ be a transitive action of $F_2(N)$ on $A$.
For each $a\in A$ and $i,j\in N$
with $i\not=j$ let 
$$
m_{i,j;a}=m_{j,i;a}\in (\ndN+1)\cup\{\infty\}
$$ be such that $\theta (i,j;a)\in\ndN\Longrightarrow{\frac {m_{i,j;a}} {\theta (i,j;a)}}\in\ndN\cup\{\infty\}$
or
$\theta (i,j;a)=\infty\Longrightarrow m_{i,j;a}=\infty$.
Set $$
\vecm:=(m_{i,j;a}\,|\,i,j\in N,i\ne j,a\in A).
$$
Let  
\begin{align}\label{eq:defofW}
W\quad =\quad (W,N,A,\sa ,\vecm )
\end{align}  be the
semigroup generated by the set
$$
\{0,e_a,s_{i,a}\,|\,a\in A ,i\in N\}
$$
and defined by the relations
\begin{gather}
  \label{eq:erel0}
\quad\quad 00=e_a0=0e_a=s_{i,a}0=0s_{i,a}=0,
\\
\label{eq:erel}
\begin{aligned}
e_a^2=&e_a, & e_ae_b=&0 \text{ for }a\not=b,\\
e_{\sig{i}{a}}s_{i,a}=&s_{i,a}e_a=s_{i,a}, & s_{i,\sig{i}{a}}s_{i,a}=&e_a,
\end{aligned}\\
\label{eq:coxrel}
\begin{aligned}
s_is_j\cdots s_js_{i,a}=s_js_i\cdots s_is_{j,a} 
 \text{ ($m_{i,j;a}$ factors)}
 & \qquad \text{if $m_{i,j;a}$ is finite and odd,}\\
s_js_i\cdots s_js_{i,a}=s_is_j\cdots s_is_{j,a}
 \text{ ($m_{i,j;a}$ factors)}
 & \qquad \text{if $m_{i,j;a}$ is finite and even,}
\end{aligned}
\end{gather} 
where we use
the convention: 
\begin{align}\label{eq:cvofgren}
s_js_{i,a}:=s_{j,i\sa a}s_{i,a},\quad
s_is_js_{i,a}:=s_{i,j\sa i\sa a}s_js_{i,a},\ldots
\end{align} See also \eqref{eq:coxgrconv} below.
\end{defin}
\subsection{Sign representation}
Let $\ndZ A$ be the free $\ndZ$-module generated by $A$,
that is, $$
\ndZ A=\oplus_{a\in A}\ndZ a.
$$
Then there exists a unique semigroup homomorphism
$$
{\widetilde {{\rm sgn}}}:W
\rightarrow{\rm{End}}_\ndZ (\ndZ A)
$$
such that
\begin{align}
{\widetilde {{\rm sgn}}}(0)(b)=&0, &
{\widetilde {{\rm sgn}}}(e_a)(b)=&\delta_{ab}b,&
{\widetilde {{\rm sgn}}}(s_{i,a})(b)=&(-1)\delta_{ab}i\sa a
\end{align}
for $a,b\in A$ and $i\in N$,
where $\delta $ means Kronecker's symbol. Hence for $w\in W$
one has $$w\ne0$$ if and only if 
$w=e_a$ for some $a\in A$ or
there exist $m\in\ndN$ and 
$i_j\in N$, $b_j\in A$ with $1\le j\le m$
such that $b_j= i_{j+1}\sa b_{j+1}$ and
$w=s _{i_1,b_1}\cdots s_{i_{m-1},b_{m-1}}s_{i_m,b_m}
$. If this is the case, we use the convention
\begin{align}\label{eq:coxgrconv}
s _{i_1}\cdots s_{i_{m-1}}s_{i_m,b_m}:=w,
\end{align} and, if $m=0$, $s _{i_1}\cdots s_{i_{m-1}}s_{i_m,a}$ means $e_a$.
We note again 
\begin{lemma}\label{lemma:bcofgr} 
{\rm{(1)}} $s _{i_1}\cdots s_{i_{m-1}}s_{i_m,a}\ne 0$ for all $a\in A$ and $m\in\ndN\cup\{0\}$. 
\par {\rm{(2)}} If $s _{i_1}\cdots s_{i_{m-1}}s_{i_m,a}=s _{j_1}\cdots s_{j_{r-1}}s_{j_r,b}$,
then $a=b$, $i_1\cdots i_m\sa a=j_1\cdots j_r\sa b$
and $(-1)^m=(-1)^r$. 
\end{lemma}

\subsection{Generalization of Root systems}

\begin{defin}\cite[Definition~2]{HY} \label{df:srs}
We call a quadruple $(R,N,\is ,\sa )$
a {\it{multi-domains root system}}
if the following conditions hold.
\begin{enumerate}
  \item \label{en:transact}
    $N$ and $\is $ are non-empty
    sets and $\sa $ is a transitive action
    of $F_2(N)$ on $\is $.
  \item \label{en:Rpi}
    Let $V_0$ be the $|N|$-dimensional $\ndR $-linear space.
    Then
    $$R =\{(R_a,\pi _a,S_a)\,|a\in \is \},$$
    where
    $\pi _a=\{\alpha _{n,a}\,|\,n\in N\}\subset R_a\subset V_0$,
    and $\pi _a$ is a basis of $V_0$ for all $a\in \is $.
  \item \label{en:R+}
    $\rs _a=\rs ^+_a\cup -\rs ^+_a$ for all $a\in \is $, where
    $\rs ^+_a=(\ndN\cup\{0\})\pi _a\cap \rs _a$.
  \item \label{en:rbeta}
    For any $i\in N$ and $a\in \is $ one has
    $\ndR \alpha _{i,a}\cap R_a=\{\alpha _{i,a},-\alpha _{i,a}\}$.
  \item \label{en:Sa}
    $S_a=\{\s_{i,a}\,|\,i\in N\}$, and for each $a\in \is $ and
    $i\in N$ one has
    $\s_{i,a}\in \GL(V_0)$,
    \begin{align*}
    \s_{i,a}(\rs _a)=&\rs _{\sig{i}{a}},&
    \s_{i,a}(\alpha _{i,a})=&-\alpha _{i,\sig{i}{a}},&
    \s_{i,a}(\alpha _{j,a})\in
    \alpha _{j,\sig{i}{a}}+(\ndN\cup\{0\})\alpha _{i,\sig{i}{a}}
    \end{align*}
    for all $j\in N\setminus \{i\}$.
  \item \label{en:tautau}
    $\s_{i,\sig{i}{a}}\s_{i,a}=\id $ for $a\in \is $ and $i\in N$.
  \item \label{en:taueq}
    Let $a\in \is $, $i,j\in N$, $i\not=j$, and
    $d=|((\ndN\cup\{0\})\alpha _{i,a}+(\ndN\cup\{0\})\alpha _{j,a})\cap \rs _a|$.
    If $d$ is finite then $\theta (i,j;a)$ is finite and it divides $d$.
\end{enumerate}
\end{defin}

{\it{Convention.}} We write 
$$
(R,N,\is ,\sa )\in\srs 
$$ if $(R,N,\is ,\sa )$ is a multi-domains root system,
that is, $\srs
=\{(R,N,\is ,\sa )\}$ denotes the family of all
the multi-domains root systems.

\begin{defin}{\rm{\cite[Definition~4]{HY}}}\label{defin:assCox} Let $(R,N,\is ,\sa )\in \srs $.
Let $\vecm:=(m_{i,j;a}\,|\,i,j\in N,i\not=j,a\in \is )$
be 
such that
$m_{i,j;a}:=|((\ndN\cup\{0\})\alpha _{i,a}+(\ndN\cup\{0\})\alpha _{j,a})\cap \rs _a|$.
Then we call $(W,N,\is ,\sa ,\vecm )$
the {\it{Coxeter groupoid
associated with $(R,N,\is ,\sa )$}}.
\end{defin}

\begin{theorem}{\rm{\cite[Theorem~1]{HY}}}\label{theorem:rep}
Let $(R,N,\is ,\sa )\in \srs $.
Set $V=\bigoplus _{a\in \is }V_a$, where
$V_a=V_0$. Let
$P_a:V\to V_a$ and $\iota _a:V_a\to V$
be the canonical projection
and the canonical inclusion map respectively. Then the assignment
$\rho :$ $0\mapsto 0\cdot\id _V$, $e_a\mapsto \iota _aP_a$,
$s_{i,a}\mapsto \iota _{\sig{i}{a}}\s _{i,a}P_a$,
gives 
a faithful representation $(\rho ,V)$
of the Coxeter groupoid
$(W,N,\is ,\sa ,\vecm )$
associated with $(R,N,\is ,\sa )$.
\end{theorem}
\subsection{Matsumoto-type theorem}

Define $\ell :W\to \ndN\cup\{0\}\cup \{-\infty \}$ to be the map
such that $\ell (0)=-\infty $, $\ell (e_a)=0$ for all $a\in \is $,
and 
\begin{align*}
 \ell (w)=\min \{m\in \ndN \,|\,&w=s_{i_1}\cdots s_{i_{m-1}}s_{i_m,a}
 \text{ for some $i_1,\ldots ,i_m\in N$, $a\in \is $}\}
\end{align*}
for all $w\in W\setminus(\{0\}\cup\{e_a|a\in A\})$; we also refer
to Lemma~\ref{lemma:bcofgr}~(1) for this definition of $\ell$.
One has
\begin{align}\label{eq:lwinv}
 \ell (w)=&\ell (w^{-1})
\end{align} for $w\in W\setminus\{0\}$, and
\begin{align}\label{eq:lww'}
 \ell (ww')\le &\ell (w)+\ell (w')
\end{align}
for $w,w'\in W$.
We say that a product $w=s_{i_1}\cdots s_{i_{m-1}}s_{i_m,a}\in W$
is 
{\it{reduced}} if $m=\ell (w)$.

\begin{defin}\cite[Definition~5]{HY}\label{defin:Wmon}
Let $W=(W,N,\is ,\sa ,\vecm )$ be a Coxeter groupoid.
Let $
\Wmon = (\Wmon ,N,\is ,\sa ,\vecm )$ denote the semigroup generated by the set
$\{0,\te_a,\ts_{i,a}\,|\,a\in \is ,i\in N\}$
and defined by the relations
\begin{gather}
00=0,\quad 0\te_a=\te_a0=0\ts_{i,a}=\ts_{i,a}0=0,\\
\te_a^2=\te_a,\quad  \te_a\te_b=0 \text{ for }a\not=b,\quad
\te_{\sig{i}{a}}\ts_{i,a}=\ts_{i,a}\te_a=\ts_{i,a},\\
\begin{aligned}
\ts_i\ts_j\cdots \ts_j\ts_{i,a}=\ts_j\ts_i\cdots \ts_i\ts_{j,a} 
 \text{ ($m_{i,j;a}$ factors)}
 & \qquad \text{if $m_{i,j;a}$ is finite and odd,}\\
\ts_j\ts_i\cdots \ts_j\ts_{i,a}=\ts_i\ts_j\cdots \ts_i\ts_{j,a}
 \text{ ($m_{i,j;a}$ factors)}
 & \qquad \text{if $m_{i,j;a}$ is finite and even.}
\end{aligned}
\end{gather}

\end{defin}

\begin{theorem} {\rm{\cite[Theorem~5]{HY}}}
{\rm{(Matsumoto-type theorem of the Coxeter groupoids)}}
\label{theorem:matsu}
Let $W=(W,N,\is ,\sa ,\vecm )$
 be the Coxeter groupoid 
 associated with $(R,N,\is ,\sa )\in \srs$
 (see Definition~\ref{defin:assCox}).
 Suppose that $m\in \ndN\cup\{0\}$, $a\in \is $,
 and $(i_1,\ldots ,i_m),(j_1,\ldots ,j_m)\in N^m$ such that
$\ell (s_{i_1}\cdots s_{i_{m-1}}s_{i_m,a})=m$ and equation
 $$s_{i_1}\cdots s_{i_{m-1}}s_{i_m,a}=
 s_{j_1}\cdots s_{j_{m-1}}s_{j_m,a}
$$ 
holds in $W$. Then in the semigroup 
$(\Wmon ,N,\is ,\sa ,\vecm )$
one has
$$ \ts_{i_1}\cdots \ts_{i_{m-1}}\ts_{i_m,a}=
 \ts_{j_1}\cdots \ts_{j_{m-1}}\ts_{j_m,a}.
$$
\end{theorem}
\begin{corollary} {\rm{\cite[Corollary~6]{HY}}}
  \label{corollary:matsu}
Let $W=(W,N,\is ,\sa ,\vecm )$
 be the Coxeter groupoid 
 associated with $(R,N,\is ,\sa )\in \srs$
 (see Definition~\ref{defin:assCox}).
 Suppose that $m\in \ndN\cup\{0\}$, $a\in \is $,
 and $(i_1,\ldots ,i_m)\in N^m$ such that 
 $\ell (s_{i_1}\cdots s_{i_{m-1}}s_{i_m,a})<m$
 holds in $(W,N,\is ,\sa ,\vecm )$.
 Then there exist $j_1,\ldots ,j_m\in N$ and $t\in \{1,\ldots ,m-1\}$
 such that $j_t=j_{t+1}$ and in the semigroup
 $(\Wmon ,N,\is ,\sa ,\vecm )$ one has the equation
 \begin{align*}
\ts_{i_1}\cdots \ts_{i_{m-1}}\ts_{i_m,a} 
= 
 \ts_{j_1}\cdots\ts_{j_t}\ts_{j_{t+1}}\cdots \ts_{j_{m-1}}\ts_{j_m,a}.
 \end{align*}
 \end{corollary}
 
In the next section, we also need
  \begin{proposition}{\rm{\cite[Corollary~3]{HY}}}
\label{proposition:l=m-1}
Let $m\in \ndN $, $(i_1,\ldots ,i_m,j)\in N^{m+1}$,
and $a\in \is $, and suppose that
$\ell (s_{i_1}\cdots s_{i_{m-1}}s_{i_m,a})=m$. 
Then:
\begin{itemize}
\item[{\rm{(1)}}] $m
 =|\s _{i_1}\cdots \s_{i_{m-1}}\s_{i_m,a}(R^+_a)
 \cap -R^+_{i_1\cdots i_m \sa a}|$.
\item[{\rm{(2)}}] $\ell (s_{i_1}\cdots s_{i_m}s_{j,j\sa a})=m-1$ 
$\Longleftrightarrow$
$\s_{i_1}\cdots \s_{i_{m-1}}\s_{i_m,a}(\alpha _{j,a})\in
-R^+_{i_1\cdots i_m\sa a}$.
\item[{\rm{(3)}}] $\ell (s_{i_1}\cdots s_{i_m}s_{j,j\sa a})=m+1$ $\Longleftrightarrow$
$\s_{i_1}\cdots \s_{i_{m-1}}\s_{i_m,a}(\alpha _{j,a})\in
R^+_{i_1\cdots i_m\sa a}$.
\end{itemize}
\end{proposition}

\begin{example}
Here we treat the finite dimensional simple Lie superalgebra 
$D(2,1;x)$, where $x\notin\{0,-1\}$. 
Note that it has 14 (positive
and negative) roots.
One has $m_{i,j;a(k)}
=2+\delta_{k0}+(1-\delta_{k0})(\delta_{ik}+\delta_{jk})$
and 
$\s_{i,a(k)}(\alpha_{j,a(k)})
=\alpha_{j,i\sa a(k)}+
\delta_{3,m_{i,j;a(k)}}\alpha_{i,i\sa a(k)}$
for $i\ne j$.
Moreover
$$
\begin{array}{ll}
  R_{a(k)}^+=&\pi_{a(k)}\cup\{\alpha_{i,a(k)}+\alpha_{j,a(k)}|m_{i,j;a(k)}=3\} \\
  &\cup\{\alpha_{1,a(k)}+\alpha_{2,a(k)}+\alpha _{3,a(k)}\} \\
  &\cup\{\alpha_{i,a(k)}+2\alpha_{k,a(k)}+\alpha_{j,a(k)}|m_{i,j;a(k)}=2\}.
\end{array}
$$ Note that $D(2,1;1)=D(2,1)=\frak{osp}(4|2)$ (see also Section~\ref{subsection:twoone}).
\begin{figure}
  \begin{center}
\setlength{\unitlength}{1mm}
\begin{picture}(80,65)(10,-8)

\put(-7,50){$a(2)$}\put(53,50){$a(1)$}
\put(8,10){$a(0)$}\put(53,10){$a(3)$}

\put(17,20){\line(0,1){15}}\put(14,27){$2$}
\put(39.5,5){\line(1,0){15}}\put(46,6.5){$3$}
\put(39.5,25){\line(1,1){15}}\put(45,35){$1$}

\put(8.5,46.5){$1$}\put(23.5,46.5){$x$}
\put(1, 40){$1$}\put(16, 40){$2$}\put(31,40){$3$}
\put(2,45){\circle{3}}
\put(17,45){\circle{3}}\put(15.70,44.20){$\times$}
\put(32,45){\circle{3}}
\put(3.5, 45){\line(1,0){12}}\put(18.5,45){\line(1,0){12}}

\put(18.5,13.5){$-1$}\put(22,1.5){$-x$}\put(33.5,11.5){$x+1$}
\put(31, 23){$1$}\put(16, 0){$2$}\put(31, 0){$3$}
\put(32,20){\circle{3}}\put(30.70,19.20){$\times$}
\put(17,5){\circle{3}}\put(15.70,4.20){$\times$}
\put(32,5){\circle{3}}\put(30.70,4.20){$\times$}
\put(18, 6){\line(1,1){13}}\put(18.5, 5){\line(1,0){12}}
\put(32, 6.5){\line(0,1){12}}

\put(68.5,46.5){$1$}\put(79,46.5){$-x-1$}
\put(61, 40){$2$}\put(76, 40){$1$}\put(91,40){$3$}
\put(62,45){\circle{3}}
\put(77,45){\circle{3}}\put(75.70,44.20){$\times$}
\put(92,45){\circle{3}}
\put(63.5, 45){\line(1,0){12}}\put(78.5,45){\line(1,0){12}}

\put(68.5,6.5){$x$}\put(79,6.5){$-x-1$}
\put(61, 0){$2$}\put(76, 0){$3$}\put(91,0){$1$}
\put(62,5){\circle{3}}
\put(77,5){\circle{3}}\put(75.70,4.20){$\times$}
\put(92,5){\circle{3}}
\put(63.5, 5){\line(1,0){12}}\put(78.5,5){\line(1,0){12}}

\end{picture}
  \end{center}
  \caption{Dynkin diagrams of the Lie superalgebra $D(2,1;x)$}
  \label{fig:D21}
\end{figure}
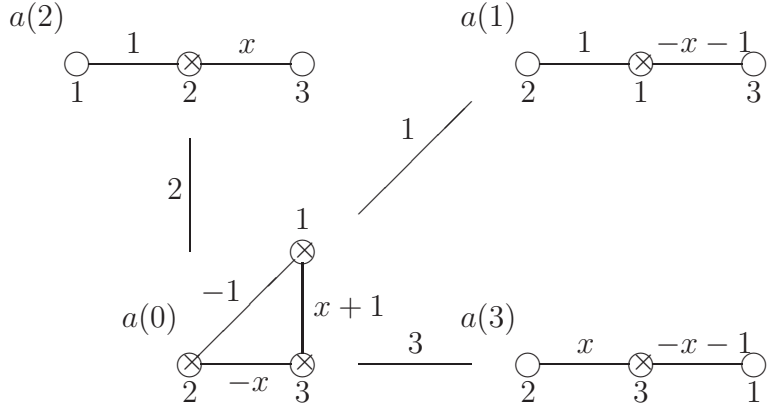
Let $w_{a(2)}:=s_{3,a(2)}s_{2,a(0)}s_{3,a(3)}s_{1,a(3)}
s_{3,a(0)}s_{2,a(2)}s_{1,a(2)}$.
Then
$\rho(w_{a(2)})=-\id _{V_{a(2)}}$.
Indeed:
\begin{align*}
&\alpha_{1,a(2)}\mapsto -\alpha_{1,a(2)}\mapsto
-\alpha_{1,a(0)}-\alpha_{2,a(0)}\mapsto
-\alpha_{1,a(3)}-\alpha_{2,a(3)}-2\alpha_{3,a(3)}\\ \nopagebreak
&\quad \mapsto -\alpha_{1,a(3)}-\alpha_{2,a(3)}-2\alpha_{3,a(3)}
\mapsto -\alpha_{1,a(0)}-\alpha_{2,a(0)}\mapsto
-\alpha_{1,a(2)}\mapsto -\alpha_{1,a(2)},\\
&\alpha_{2,a(2)}\mapsto\alpha_{1,a(2)}+\alpha_{2,a(2)}\mapsto\alpha_{1,a(0)}
\mapsto\alpha_{1,a(3)}+\alpha_{3,a(3)}\mapsto\alpha_{3,a(3)}
\mapsto -\alpha_{3,a(0)}\\ \nopagebreak
&\quad \mapsto -\alpha_{2,a(2)}-\alpha_{3,a(2)}\mapsto -\alpha_{2,a(2)},\\
&\alpha_{3,a(2)}\mapsto\alpha_{3,a(2)}\mapsto\alpha_{2,a(0)}+\alpha_{3,a(0)}
\mapsto\alpha_{2,a(3)}\mapsto\alpha_{2,a(3)}
\mapsto\alpha_{2,a(0)}+\alpha_{3,a(0)}\\ \nopagebreak
&\quad \mapsto\alpha_{3,a(2)}\mapsto -\alpha_{3,a(2)}.
\end{align*}
By Proposition~\ref{proposition:l=m-1}(1),
we have $\ell(w_{a(2)})=|R_{a(2)}^+|=7$,
$w_{a(2)}$ is the longest word.
Let
$w':=(s_{3,a(2)}s_{2,a(0)})^{-1}w_{a(2)}$.
Then $\ell (w')=5$.
By Theorem~\ref{theorem:matsu},
$w'$
has the following 
four reduced expressions:
\begin{align*}
w'=&{\underline{s_{3,a(3)}s_{1,a(3)}
s_{3,a(0)}}}s_{2,a(2)}s_{1,a(2)}
=s_{1,a(1)}s_{3,a(1)}{\underline{s_{1,a(0)}s_{2,a(2)}s_{1,a(2)}}}\\
=&s_{1,a(1)}{\underline{s_{3,a(1)}s_{2,a(1)}}}s_{1,a(0)}s_{2,a(2)}
=s_{1,a(1)}s_{2,a(1)}s_{3,a(1)}s_{1,a(0)}s_{2,a(2)}.
\end{align*}
\end{example}

\section{Main theorems---Irreducible representations
of the Iwahori-Hecke type algeras 
$H_q(A(m,n))$, $H_q(B(m,n))$, $H_q(C(n))$
and $H_q(D(m,n))$
associated with the 
Lie superalgebras $A(m,n)$, $B(m,n)$, $C(n)$, $D(m,n)$}

\subsection{Definition of Lie superalgebras}
As for the terminology concerning
Lie superalgebras, we refer to
\cite{kac1}.

Let $\mathfrak{v}=\mathfrak{v}(0)\oplus\mathfrak{v}(1)$
be a $\ndZ /2\ndZ $-graded
$\ndC $-linear space. If $i\in \{0,1\}$ and $j\in \ndZ $ such that
$j-i\in 2\ndZ $ then let $\mathfrak{v}(j)=\mathfrak{v}(i)$.
If $X\in\mathfrak{v}(0)$ (resp.\ $X\in\mathfrak{v}(1)$)
then we write 
\begin{equation}\label{eqn:pxzeroone}
\mbox{$\dg (X)=0$ (resp.\ $\dg (X)=1$)}
\end{equation}
and
we say that $X$ is an {\it{even}} (resp.\ {\it{odd}})
element.
If $X\in\mathfrak{v}(0)\cup\mathfrak{v}(1)$, then we say
that $X$ is a
{\it {homogeneous}} element and that
$\dg (X)$ is the {\it{parity}} (or {\it{degree}}) of $X$.
If $\mathfrak{w}\subset \mathfrak{v}$ is a subspace and
$\mathfrak{w}=(\mathfrak{w}\cap\mathfrak{v}(0))\oplus
(\mathfrak{w}\cap\mathfrak{v}(1))$
(resp. $\mathfrak{w}\subset\mathfrak{v}(0)$,
resp. $\mathfrak{w}\subset\mathfrak{v}(1)$), then we say that
$\mathfrak{w}$ is a {\it{graded}}
(resp. {\it{even}}, resp. {\it{odd}})
subspace.

Let $\mathfrak{g}=\mathfrak{g}(0)\oplus\mathfrak{g}(1)$ be a
$\ndZ /2\ndZ $-graded
$\ndC $-linear space equipped with a bilinear
map $[\,,\,]:\mathfrak{g}\times\mathfrak{g}\rightarrow\mathfrak{g} $
such that $[\mathfrak{g}(i),\mathfrak{g}(j)]\subset \mathfrak{g}(i+j)$
($i$, $j\in\ndZ $); we
recall from the above paragraph that
\begin{equation}\label{eqn:defohai}
\mathfrak{g}(i)=\{X\in\mathfrak{g}\,|\,\dg (X)=i\}. 
\end{equation}
We say that $\mathfrak{g}=(\mathfrak{g},[\,,\,])$ is a
($\ndC $-){\it{Lie superalgebra}} if
for all homogeneous elements $X$, $Y$, $Z$ of $\mathfrak{g}$ the following equations
hold.
\begin{align*}
	[Y,X]=&\,-(-1)^{\dg (X)\dg (Y)}[X,Y],& &\text{(skew-symmetry)} \\
	[X,[Y,Z]]=&\,[[X,Y],Z]+(-1)^{\dg (X)\dg (Y)}[Y,[X,Z]]. & &
	\text{(Jacobi identity)}
\end{align*} We call the Lie algebra $\mathfrak{g}(0)$
the {\it {even part}} of $\mathfrak{g}$.

\subsection{
Lie superalgebras $\frak{gl}(m+1|n+1)$ and $\frak{osp}(m|n)$}
\label{subsection:twoone}

Let $m$, $n\in\ndN\cup\{0\}$.
Let:
$$
\mathcal{D}_{m+1|n+1}:=\{(p_1,\ldots,p_{m+n+2})\in\ndZ ^{m+n+2}\,|\,p_i\in\{0,1\},
\sum_{i=1}^{m+n+2}p_i=n+1\}.
$$ 
For $i,j\in\{1,\ldots,m+n+2\}$, let ${\bf{E}}_{i,j}$ denote the $(m+n+2)\times (m+n+2)$ matrix
having 1 in $(i,j)$ position and 0 otherwise, that is,
the $(i,j)$-matrix unit. 
Let ${\bf{E}}_{m+n+2}$ denote the $(m+n+2)\times (m+n+2)$ unit matrix, that is,
${\bf{E}}_{m+n+2}=\sum_{i=1}^{m+n+2}{\bf{E}}_{i,i}$.
Denote by ${\rm{M}}_{m+n+2}(\ndC )$ 
the $\ndC$-linear space of the $(m+n+2)\times (m+n+2)$-matrices,
i.e., ${\rm{M}}_{m+n+2}(\ndC )=\oplus_{i,j=1}^{m+n+2}\ndC {\bf{E}}_{i,j}$.

Let $d=(p_1,\ldots,p_{m+n+2})\in \mathcal{D}_{m+1|n+1}$.
The  
Lie superalgebra $\mathfrak{gl}(m+1|n+1)=\frak{gl}(d)$ is defined by
$\frak{gl}(d)={\rm{M}}_{m+n+2}(\ndC )$
(as a $\ndC$-linear space),
\begin{align}\label{eq:evoodptofgl}
\frak{gl}(d)(0)=\oplus_{1\leq p_i=p_j\leq m+n+2}\ndC {\bf{E}}_{i,j}, \quad
\frak{gl}(d)(1)=\oplus_{1\leq p_i\ne p_j\leq m+n+2}\ndC {\bf{E}}_{i,j}, 
\end{align}
and $[X,Y]=XY-(-1)^{r_1r_2}YX$
for $X\in\mathfrak{gl}(d)(r_1)$
and $Y\in\mathfrak{gl}(d)(r_2)$,
\newline\newline
where $XY$ and $YX$ mean the matrix product, that is,
${\bf{E}}_{i,j}{\bf{E}}_{k,l}=\delta_{j,k}{\bf{E}}_{i,l}$.
Define the $\ndC$-linear map ${\rm{str}}:\mathfrak{gl}(d)\to\ndC$
by ${\rm{str}}({\bf{E}}_{i,j})=\delta_{i,j}
(-1)^{p_i}$.
The Lie subsuperalgebra $\{X\in \mathfrak{gl}(d)\,|\,{\rm{str}}(X)=0\}$
of $\mathfrak{gl}(d)$ is denoted as $\mathfrak{sl}(m+1|n+1)
=\mathfrak{sl}(d)$.
The finite dimensional simple Lie superalgebra $A(m,n)$ 
is defined as follows.
Let $\mathfrak{z}$ be the one dimensional ideal 
$\ndC {\bf{E}}_{m+n+2}$ of $\mathfrak{gl}(d)$.
If $m\ne n$, then $A(m,n)$ means $\mathfrak{sl}(d)$.
On the other hand, $A(n,n)$ means $\mathfrak{sl}(d)/\mathfrak{z}$, and
is also denoted as
$\psl(n+1|n+1)$. 

Let $d=(p_1,\ldots,p_{m+2n})\in \mathcal{D}_{m|2n}$. Define the map
$\theta:\{1,\ldots,m+2n\}\rightarrow\{1,\ldots,m+2n\}$
by $\theta(i)=m+2n+1-i$. 
Assume that $p_{\theta(i)}=p_i$.
Let $g_i\in\{1,-1\}$ be such that 
$g_i=-1$ if $p_i=1$ and $i<\theta(i)$ and 
$g_i=1$ otherwise. We have an automorphism
$\Omega$
of $\mathfrak{gl}(d)$ defined by
$\Omega ({\bf{E}}_{i,j})=-(-1)^{p_ip_j+p_j}g_ig_j{\bf{E}}_{\theta(j),\theta(i)}$.
The Lie superalgebra $\frak{osp}(m|2n)$ means 
$\{X\in\mathfrak{gl}(d)|\Omega (X)=X\}$. We also denote
$\frak{osp}(m|2n)$ as follows:
 \begin{align*}
B(m-1,n)
&= 
\frak{osp}(2m-1|2n) &\quad\text{if $m$, $n\in\ndN$,}
\\
D(m+1,n)
&= 
\frak{osp}(2m+2|2n) &\quad\text{if $m$, $n\in\ndN$,}
\\
C(n+1) &=
\frak{osp}(2|2n) &\quad\text{if $n\in\ndN$.}
 \end{align*}
We also note that $\frak{osp}(2m+1|0)$,
$\frak{osp}(0|2n)$, and $\frak{osp}(2m|0)$ are isomorphic to the
simple Lie algebras of type $B_m$ (if $m\geq 2$),  $C_n$ (if $n\geq 3$) and 
$D_m$  (if $m\geq 4$)
respectively, so $\frak{osp}(2m+1|0)\cong\frak{o}_{2m+1}$, 
$\frak{osp}(0|2n)\cong\frak{sp}_{2n}$ and $\frak{osp}(2m|0)\cong\frak{o}_{2m}$.
As for the even part $\frak{osp}(m|2n)(0)$ of $\frak{osp}(m|2n)$, we have
\begin{align}\label{eq:evenpartofosp}
\frak{osp}(m|2n)(0)\cong \frak{osp}(m|0)\oplus \frak{osp}(0|2n).
\end{align}

\subsection{Definition of Iwahori-Hecke type algebras}

\begin{defin}\label{df:defIHalgebroid}
Let $W = (W,N,A,\sa ,\vecm )$ be the groupoid  introduced in \eqref{eq:defofW}.
Assume that $A$ is finite. 
Let $q\in\ndC$.
Let $H_q(W)$ be the $\ndC$-algebra (with $1$) generated by 
\begin{align}\label{eq:genHalgebroid}
\{E_a,T_{i,a}|a\in A,i\in N\}
\end{align} and defined by the relations
\begin{align}
E_a^2=E_a, & &  \label{eq:HIrel1-1} \\
E_{i\sa a}T_{i,a}E_a=T_{i,a}, & &  \label{eq:HIrel1-2} \\ 
 \sum_{a\in A}E_a=1 & &  \label{eq:HIrel1-3} \\
E_aE_b=0 & & \mbox{if $a\ne b$}, & \label{eq:HIrel2}  \\
(T_{i,a}-qE_a)(T_{i,a}+E_a)=0 & &  \mbox{if $i\sa a=a$}, & \label{eq:HIrel3} \\
T_{i,i\sa a}T_{i,a}=E_a & & \mbox{if $i\sa a\ne a$}, & \label{eq:HIrel4}  \\
T_iT_j\cdots T_jT_{i,a}=T_jT_i\cdots T_iT_{j,a} 
 \text{ ($m_{i,j;a}$ factors)}
 & & \text{if $m_{i,j;a}$ is finite and odd,} &\label{eq:HIrel5}  \\
T_jT_i\cdots T_jT_{i,a}=T_iT_j\cdots T_iT_{j,a}
 \text{ ($m_{i,j;a}$ factors)}
 & & \text{if $m_{i,j;a}$ is finite and even,} & \label{eq:HIrel6} 
\end{align} where, in \eqref{eq:HIrel5}-\eqref{eq:HIrel6}, we use the same convention as that of \eqref{eq:cvofgren}
with $s_{i,a}$ in place of $T_{i,a}$.
\end{defin}

\begin{lemma} \label{lemma:spanHq}
Let $W=(W,N,\is ,\sa ,\vecm )$ be the Coxeter groupoid associated with an element $(R,N,\is ,\sa )$
of $\srs $ (see Definition~\ref{defin:assCox}).
Assume that $A$ is finite.
Then there exists a map $f:W\rightarrow H_q(W)$ such that 
\begin{align}
f(0)=0, \,\, f(e_a)=E_a, & & & \label{eq:spanHq1} \\
f(s_{i,a}w) =T_{i,a}f(w) & & \mbox{if $w\in W\setminus\{0\}$ and $\ell ( s_{i,a}w)=1+\ell (w)$}. &
\label{eq:spanHq2}
\end{align} Further, as a $\ndC$-linear space, $H_q(W)$
is spanned by $f(W\setminus\{0\})$. In particular,
if $W$ is finite, then
\begin{align}\label{eq:spanHq3}
\dim H_q(W)\leq |W|-1.
\end{align}
\end{lemma}

{\it{Proof.}} Let $\Wmon$ be the semigroup introduced in Definition~\ref{defin:Wmon}
for $W$.
It is easy to show that there exists a unique semigroup homomorphism
${\widetilde{f}}:{\widetilde{W}}\rightarrow H_q(W)$ such that 
${\widetilde{f}}(0)=0$, ${\widetilde{f}}(\te_a)=E_a$
and ${\widetilde{f}}(\ts_{i,a})=T_{i,a}$.
By Theorem~\ref{theorem:matsu}, there exists a unique map
$f:W\rightarrow H_q(W)$ such that $f(0)=0$
and $f(w)={\widetilde{f}}(\ts_{i_1}\cdots \ts_{i_{m-1}}\ts_{i_m,a})$ if $w\in W\setminus\{0\}$,
$\ell(w)=m$ and $w=s_{i_1}\cdots s_{i_{m-1}}s_{i_m,a}$.
Then $f$ satisfies \eqref{eq:spanHq1}-\eqref{eq:spanHq2}, as desired.

We show 
\begin{align}\label{eq:spanHq4}
\forall w \in W, \forall i \in N, \forall a \in A,  T_{i,a}f(w)\in \ndC f(s_{i,a}w)+\ndC f(w).
\end{align}  If $s_{i,a}w=0$, then clearly $T_{i,a}f(w)=0$ holds.
If $w\ne 0$, $s_{i,a}w\ne 0$ and $\ell ( s_{i,a}w)= 1+\ell (w)$,
then \eqref{eq:spanHq4} follows from \eqref{eq:spanHq2}.
Assume that $w\ne 0$, $s_{i,a}w\ne 0$ and $\ell ( s_{i,a}w)\ne 1+\ell (w)$.
Then by \eqref{eq:lwinv} and Proposition~\ref{proposition:l=m-1}, we have
$\ell ( s_{i,a}w)= \ell (w)-1$, so $f(w)=T_{i,i\sa a}f(s_{i,a}w)$.
Since $T_{i,a}f(w)=T_{i,a}T_{i,i\sa a}f(s_{i,a}w)$,
we have $T_{i,a}f(w)=f(s_{i,a}w)$ if $i\sa a \ne a$,
and $T_{i,a}f(w)=(q-1)f(w)+qf(s_{i,a}w)$ otherwise. 
Hence we have \eqref{eq:spanHq4}, as desired.

It is clear from \eqref{eq:spanHq4} that the rest of the statement follows.
\hfill $\Box$

\begin{notation}
Let $r\in\ndN$. Let $V^{(r)}_0$ be the $r$-dimensional $\ndR$-linear space
with a basis $\{\varepsilon_i|1\leq i\leq r\}$.
Let $V^{(r),\prime}_0$ be the subspace of $V^{(r)}_0$ formed 
by the elements $\sum_{i=1}^rx_i\varepsilon_i$ with
$x_i\in\ndR$ and $\sum_{i=1}^rx_i=0$, so $\dim V^{(r),\prime}_0=r-1$.
For a non-zero element $x=\sum_{i=1}^{|N|}x_i\varepsilon_i$ of $V^{(r)}_0$
with $x_i\in\ndR$, define ${\widetilde{\sigma}}_x\in{\rm{GL}} (V^{(r)}_0)$ by 
${\widetilde{\sigma}}_x (\varepsilon_j)=
\varepsilon_j-2x_j(\sum_{i=1}^{|N|}x_i^2)^{-1}x$, that is,
${\widetilde{\sigma}}_x$ is the reflection of $V^{(r)}_0$ with respect to 
the hyperplane of $V^{(r)}_0$ orthogonal to $x$. Note that
if $x\in V^{(r),\prime}_0$, then ${\widetilde{\sigma}}_x(V^{(r),\prime}_0)=V^{(r),\prime}_0$. 
\end{notation}
 
\subsection{Basic of Iwahori-Hecke algebras}
For the basic facts about the Iwahori-Hecke algebras, we refer to \cite{GU}.
Let $W=(W,N,A,\sa ,\vecm )$ be the groupoid introdued in \eqref{eq:defofW}.
In this subsection we always assume that 
\begin{align}
 \mbox{{\it{$|A|=1$ and $N$ and $W$ are finite}}}. 
\end{align} Let $a\in A$, so $A=\{a\}$. Then $W\setminus\{0\}$ is nothing but
the Coxeter group associated with the Coxeter system 
$(W\setminus\{0\},\{s_{i,a}|,i\in N\})$. 
In this case, we also denote 
$H_q(W)$ and $T_{i,a}$ by 
$H_q(W\setminus\{0\})$ and $T_i$ respectively.
That is, $H_q(W\setminus\{0\})$
is the $\ndC$-algebra
(with $1$) generated by $T_i$ ($i\in N$)
and defined by the relations
\begin{gather}
\begin{aligned}
  \label{eq:IHone}
(T_i-q)(T_i+1)=0,
\end{aligned}\\
\label{eq:IHtwo}
\begin{aligned}
T_iT_j\cdots T_jT_i=T_jT_i\cdots T_iT_j 
 \text{ ($m_{i,j;a}$ factors)}
 & \qquad \text{if $m_{i,j;a}$ is odd,}\\
T_jT_i\cdots T_jT_i=T_iT_j\cdots T_iT_j
 \text{ ($m_{i,j;a}$ factors)}
 & \qquad \text{if $m_{i,j;a}$ is even.}
\end{aligned}
\end{gather} It is well-known that $\dim H_q(W\setminus\{0\})=|W\setminus\{0\}|$.
In this paper we fix a complete set of non-equivalent irreducible representations
of $H_q(W\setminus\{0\})$ by 
\begin{align}\label{eq:IHthree}
\{\rho_{H_q(W\setminus\{0\}),\lambda}:H_q(W\setminus\{0\})\to 
\End_\ndC (V_{H_q(W\setminus\{0\}),\lambda}) |\lambda\in\Lambda_{H_q(W\setminus\{0\})}\},
\end{align}
where $\Lambda_{H_q(W\setminus\{0\})}$ is an index set. Define the
polynomial $P_{W\setminus\{0\}}(q)$ in $q$ by
\begin{align}\label{eq:IHfour}
P_{W\setminus\{0\}}(q):=\sum_{w \in W\setminus\{0\}}q^{\ell(w)}.
\end{align} This is called the {\it{Poincar\'{e} polynomial}} of $W\setminus\{0\}$.

It is well-known \cite{GU} (see also \cite[(25.22) and (27.4)]{CR}) that for $q\in\ndC\setminus\{0\}$,
the following three conditions are equivalent.
\newline\par
(i) $P_{W\setminus\{0\}}(q)\ne 0$ holds. \par
(ii) $H_q(W\setminus\{0\})$ is a semisimle algebra. \par
(iii) The map 
\begin{align}\label{eq:IHfive}
\bigoplus_{\lambda\in\Lambda_{H_q(W\setminus\{0\})}}\rho_{H_q(W\setminus\{0\}),\lambda}:H_q(W\setminus\{0\})
\rightarrow \bigoplus_{\lambda\in\Lambda_{H_q(W\setminus\{0\})}}
\End_\ndC (V_{H_q(W\setminus\{0\}),\lambda})
\end{align} defined by  $X\mapsto \oplus_{\lambda\in\Lambda_{H_q(W\setminus\{0\})}}\rho_{H_q(W\setminus\{0\}),\lambda}(X)$
is a $\ndC$-algebra isomorphism.
\newline\newline In particular, 
\begin{align}\label{eq:IHsix}
q\cdot P_{W\setminus\{0\}}(q)\ne 0 \Longrightarrow \dim
H_q(W\setminus\{0\}) =\sum_{\lambda\in\Lambda_{H_q(W\setminus\{0\})}}(\dim V_{H_q(W\setminus\{0\}),\lambda})^2.
\end{align}

Assume that $N=\{1,2,\ldots,n\}$ and $m_{i,i+1;a}=3$ and 
$m_{i,j;a}=2$ ($|j-i|\geq 2$). Then  $W$ is the Coxeter groupoid 
associated with $(R,N,\is ,\sa )\in \srs $ such that 
$V_0=V^{(n+1),\prime}_0$, $R_a^+=\{\varepsilon_i-\varepsilon_j|1\leq i<j\leq n+1\}$,
$\al_{i,a}=\varepsilon_i-\varepsilon_{i+1}$ and $\sigma_{i,a}={\widetilde{\sigma}}_{\al_{i,a}}$.
As a group, 
$W\setminus\{0\}$ is isomorphic to the symmetric group $S_{n+1}$,
so we also denote $W\setminus\{0\}$ by $S_{n+1}$ by abuse of notation. 
Note that $\dim H_q(S_{n+1})=(n+1)!$.

Assume that $N=\{1,2,\ldots,n\}$ and $m_{i,i+1;a}=3$ ($1\leq i\leq n -3$),
$m_{n-1,n;a}=4$ and 
$m_{i,j;a}=2$ ($|j-i|\geq 2$). Then  $W$ is the Coxeter groupoid 
associated with $(R,N,\is ,\sa )\in \srs $ such that 
$V_0=V^{(n)}_0$, $R_a^+=\{\varepsilon_i-\varepsilon_j, \varepsilon_i+\varepsilon_j|1\leq i<j\leq n \}
\cup\{\varepsilon_i|1\leq i\leq n\}$,
$\al_{i,a}=\varepsilon_i-\varepsilon_{i+1}$ ($1\leq i\leq n-1$), $\al_{n,a}=\varepsilon_n$
and $\sigma_{i,a}={\widetilde{\sigma}}_{\al_{i,a}}$.
We also denote $W\setminus\{0\}$ by $W(B_n)$
and $W(C_n)$. Note that $\dim H_q(W(B_n))=2^n n!$.

Assume that $N=\{1,2,\ldots,n\}$ and $m_{i,i+1;a}=3$ ($1\leq i\leq n -2$),
$m_{n-1,n;a}=2$, $m_{n-2,n;a}=3$ and 
$m_{i,j;a}=2$ ($|j-i|\geq 2$ and $1\leq i\leq n-3$).
We also denote $W\setminus\{0\}$ by $W(D_n)$. Note that $\dim H_q(W(D_n))=2^{n-1} n!$.
Then  $W$ is the Coxeter groupoid 
associated with $(R,N,\is ,\sa )\in \srs $ such that 
$V_0=V^{(n)}_0$, $R_a^+=\{\varepsilon_i-\varepsilon_j, \varepsilon_i+\varepsilon_j|1\leq i<j\leq n \}$,
$\al_{i,a}=\varepsilon_i-\varepsilon_{i+1}$ ($1\leq i\leq n-1$), $\al_{n,a}=\varepsilon_{n-1}+\varepsilon_n$
and $\sigma_{i,a}={\widetilde{\sigma}}_{\al_{i,a}}$.

It is well-known (cf. \cite[Theoerem~10.2.3 and Proposition~10.2.5]{C})
that \begin{align}\label{eq:IHseven}
P_{S_{n+1}}(q)&=\prod_{r=1}^n{\frac {q^{r+1} -1} {q-1}}, \\
P_{W(B_n)}(q)&=\prod_{r=1}^n{\frac {q^{2r} -1} {q-1}}, \label{eq:IHeight} \\
P_{W(D_n)}(q)&={\frac {q^n -1} {q-1}}\prod_{r=1}^{n-1}{\frac {q^{2r} -1} {q-1}}. \label{eq:IHnine}
\end{align}

\subsection{Iwahori-Hecke type algebra
$H_q(A(m,n))$
associated with the 
Lie superalgebra $A(m,n)$} \label{subsection:Amn}

Let 
$$
\tr : S_{m+n+2}\times \mathcal{D}_{m+1|n+1} \to \mathcal{D}_{m+1|n+1} 
$$
denote the usual (left) action of the symmetric group
$S_{m+n+2}$ on $\mathcal{D}_{m+1|n+1}$ by permutations, that is, for $\sigma\in S_{m+n+2}$,
$$
\sigma\tr(p_1,\ldots,p_{m+n+2})=
(p_{\sigma^{-1}(1)},\ldots,p_{\sigma^{-1}(m+n+2)}).
$$ 

Let $\sigma_i:=(i,i+1)\in S_{m+n+2}$. Let $W$ be the Coxeter groupoid 
associated with $(R,N,\is ,\sa )\in \srs $ such that
$N=\{1,2,\ldots,m+n+1\}$, $A=\mathcal{D}_{m+1|n+1}$, $i\sa d= \sigma_i \tr d$,
$V_0=V^{(m+n+2),\prime}_0$, $R_d^+=\{\varepsilon_i-\varepsilon_j|1\leq i<j\leq m+n+1\}$,
$\al_{i,d}=\varepsilon_i-\varepsilon_{i+1}$ and $\sigma_{i,d}={\widetilde{\sigma}}_{\al_{i,d}}$.
Denote $H_q(W)$ by $H_q(A(m,n))$. 
Then $H_q(A(m,n))$
is the $\ndC$-algebra (with $1$) generated by
\begin{align}
\{E_d\,|\,d\in\mathcal{D}_{m+1|n+1} \}
\cup\{T_{i,d}\,|\,1\leq i\leq m+n+1,\,d\in\mathcal{D}_{m+1|n+1} \}
	\label{eq:HAgen}
\end{align}
and defined by the relations \eqref{eq:HIrel1-1}-\eqref{eq:HIrel4} and
the relations
\begin{align}
	T_{i,\sigma_j\sigma_i\tr d}T_{j,\sigma_i\tr d}T_{i,d}=&
	T_{j,\sigma_i\sigma_j\tr d}T_{i,\sigma_j\tr d}T_{j,d}&  
    &\text{if $|i-j|= 1$,} \label{eq:HArel9}
	\\
	T_{i,\sigma_j\tr d}T_{j,d}=&T_{j,\sigma_i\tr d}T_{i,d}&
	&\text{if $|i-j|\geq 2$.}
	\label{eq:HArel10}
\end{align}

Define $d_e$, $d_o\in\mathcal{D}_{m+1|n+1}$ by
\begin{align}\label{eq:dedo}
d_e:=(\overbrace{0,\ldots,0}^{m+1},\overbrace{1,\ldots,1}^{n+1}),\,\,
d_o:=(\overbrace{1,\ldots,1}^{n+1},\overbrace{0,\ldots,0}^{m+1}).
\end{align}
For $d=(p_1,\ldots,p_{m+n+2})\in\mathcal{D}_{m+1|n+1}$,
define the two elements 
\begin{align} \label{eq:taudplusminus}
\tau_ {+.d}, \tau_ {-.d}\in S_{m+n+2}
\end{align} by 
\begin{align}
p_{\tau_ {\pm.d} (i)}={\frac {1\mp 1} 2} & \quad\text{and}\quad\tau_ {\pm.d} (i) \leq \tau_ {\pm.d} (j)
& &\quad\text{if $1\leq i\leq j \leq m+1$,}  \\
p_{\tau_ {\pm.d} (i)}={\frac {1\pm 1} 2} & \quad\text{and}\quad\tau_ {\pm.d} (i) \leq \tau_ {\pm.d} (j)
& &\quad\text{if $m+2\leq i\leq j \leq m+n+2$.}  
\end{align} Then $\tau_ {+,d}$ (resp. $\tau_ {-,d}$) is the minimal length element among
the elements $\sigma\in S_{m+n+2}$ satisfying the condition that
for any $i$, $i$-th component of $d_e$ (resp. $d_o$) 
is the same as $\sigma(i)$-th component $p_{\sigma(i)}$ of
$d$. 

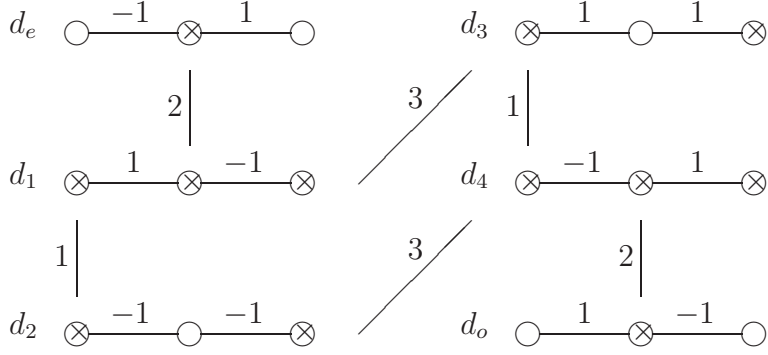
\begin{figure}
  \begin{center}
\setlength{\unitlength}{1mm}
\begin{picture}(80,50)(10,-8)

\put(-7,45){$d_e$}\put(53,45){$d_3$}
\put(-7,25){$d_1$}\put(53,25){$d_4$}
\put(-7,5){$d_2$}\put(53,5){$d_o$}

\put(17,30){\line(0,1){10}}\put(14,34){$2$}
\put(2,10){\line(0,1){10}}\put(-1,14){$1$}
\put(62,30){\line(0,1){10}}\put(59,34){$1$}
\put(77,10){\line(0,1){10}}\put(74,14){$2$}

\put(39.5,5){\line(1,1){15}}\put(46,15){$3$}
\put(39.5,25){\line(1,1){15}}\put(46,35){$3$}

\put(6.5,46.5){$-1$}\put(23.5,46.5){$1$}
\put(2,45){\circle{3}}
\put(17,45){\circle{3}}\put(15.70,44.20){$\times$}
\put(32,45){\circle{3}}
\put(3.5, 45){\line(1,0){12}}\put(18.5,45){\line(1,0){12}}

\put(8.5,26.5){$1$}\put(21.5,26.5){$-1$}
\put(2,25){\circle{3}}\put(0.70,24.20){$\times$}
\put(17,25){\circle{3}}\put(15.70,24.20){$\times$}
\put(32,25){\circle{3}}\put(30.70,24.20){$\times$}
\put(3.5, 25){\line(1,0){12}}\put(18.5,25){\line(1,0){12}}

\put(6.5,6.5){$-1$}\put(21.5,6.5){$-1$}
\put(2,5){\circle{3}}\put(0.70,4.20){$\times$}
\put(17,5){\circle{3}}
\put(32,5){\circle{3}}\put(30.70,4.20){$\times$}
\put(3.5, 5){\line(1,0){12}}\put(18.5,5){\line(1,0){12}}

\put(68.5,46.5){$1$}\put(83.5,46.5){$1$}
\put(62,45){\circle{3}}\put(60.70,44.20){$\times$}
\put(77,45){\circle{3}}
\put(92,45){\circle{3}}\put(90.70,44.20){$\times$}
\put(63.5, 45){\line(1,0){12}}\put(78.5,45){\line(1,0){12}}

\put(66.5,26.5){$-1$}\put(83.5,26.5){$1$}
\put(62,25){\circle{3}}\put(60.70,24.20){$\times$}
\put(77,25){\circle{3}}\put(75.70,24.20){$\times$}
\put(92,25){\circle{3}}\put(90.70,24.20){$\times$}
\put(63.5, 25){\line(1,0){12}}\put(78.5,25){\line(1,0){12}}

\put(68.5,6.5){$1$}\put(81.5,6.5){$-1$}
\put(62,5){\circle{3}}
\put(77,5){\circle{3}}\put(75.70,4.20){$\times$}
\put(92,5){\circle{3}}
\put(63.5, 5){\line(1,0){12}}\put(78.5,5){\line(1,0){12}}

\end{picture}
  \end{center}
  \caption{Dynkin diagrams of the Lie superalgebra $A(1,1)$}
  \label{fig:A11}
\end{figure}

\begin{example} Assume that $m=n=1$.
Then $\mathcal{D}_{2|2}=\{d_e=(0,0,1,1),d_1=(0,1,0,1),d_2=(1,0,0,1),d_3=(0,1,1,0),d_4=(1,0,1,0),d_o=(1,1,0,0)\}$.
Then 
$\tau_ {+.d_e}={{1234} \brack {1234}}$, $\tau_ {-.d_e}={{1234} \brack {3412}}$, 
$\tau_ {+.d_1}={{1234} \brack {1324}}$, $\tau_ {-.d_1}={{1234} \brack {2413}}$, 
$\tau_ {+.d_2}={{1234} \brack {2314}}$, $\tau_ {-.d_2}={{1234} \brack {1423}}$,
$\tau_ {+.d_3}={{1234} \brack {1423}}$, $\tau_ {-.d_3}={{1234} \brack {2314}}$,
$\tau_ {+.d_4}={{1234} \brack {2413}}$, $\tau_ {-.d_4}={{1234} \brack {1324}}$,
$\tau_ {+.d_o}={{1234} \brack {3412}}$, $\tau_ {-.d_o}={{1234} \brack {1234}}$. 
See also Figure~\ref{fig:A11}.
\end{example}

Now we consider $|W|$. Recall $\rho$ and $d_e$ from Theorem~\ref{theorem:rep}
and \eqref{eq:taudplusminus}
respectively. It is easy to see that
$P_{d_e}\rho(e_{d_e}We_{d_e})\iota_{d_e}\subset\{(\sum_{i=1}^{m+n+2}{\bf{E}}_{\sigma(i)i})_{|V^{(m+n+2),\prime}_0}|
\sigma\in S_{m+n+2}, \sigma(\{1,\ldots,m+1\})=\{1,\ldots,m+1\}\}$. Hence $|e_{d_e}We_{d_e}|\leq (m+1)!(n+1)!$
by Theorem~\ref{theorem:rep}, so 
$|W\setminus\{0\}|=|\mathcal{D}_{m+1|n+1}|^2|e_{d_e}W e_{d_e}|\leq {\frac {((m+n+2)!)^2} {(m+1)!(n+1)!}}$.
Hence by \eqref{eq:spanHq3}, we conclude
\begin{align}\label{eq:spanHq3typeA}
\dim H_q(A(m,n))\leq {\frac {((m+n+2)!)^2} {(m+1)!(n+1)!}}.
\end{align}

%

\begin{proposition}\label{proposition:repA}
Let $V$ and $W$ be finite dimensional $\ndC$-linear spaces, and
let ${\bf l}:H_q(S_{m+1})\rightarrow\End_\ndC (V)$ and ${\bf r}:H_q(S_{n+1})\rightarrow\End_\ndC (W)$
be $\ndC$-algebra homomorphisms, i.e., representations.
Let ${\bf l}\otimes{\bf r}:H_q(S_{m+1})\otimes H_q(S_{n+1})\rightarrow 
\End_\ndC (V\otimes W)$ denote the tensor representation of 
${\bf l}$ and ${\bf r}$ in the ordinary sense.
Let $C_{V\otimes W;d}$ be copies of 
the $\ndC$-linear space
$V\otimes W$, indexed by $d\in\mathcal{D}_{m+1|n+1} $.
Let $C_{V\otimes W}:=\oplus_{d\in\mathcal{D}_{m+1|n+1}}C_{V\otimes W;d}$. 
Let $P_d:C_{V\otimes W}\rightarrow C_{V\otimes W;d}$
and $\iota_d:C_{V\otimes W;d}\rightarrow C_{V\otimes W}$ denote 
the canonical
projection and the canonical inclusion map respectively.
Then there exists a unigue $\ndC$-algebra homomorphism 
${\bf l}\boxtimes^{A(m,n)}{\bf r}:H_q(A(m,n))\rightarrow \End_\ndC 
(C_{V\otimes W})$ satisfying the following conditions:
\newline\par
{\rm{(i)}} For each $d\in\mathcal{D}_{m+1|n+1} $,
one has $({\bf l}\boxtimes^{A(m,n)}{\bf r})(E_d)=\iota_d\circ P_d$, 
 \par
{\rm{(ii)}} For each $i\in\{1,\ldots,m+n+1\}$ and 
each $d=(p_1,\ldots,p_{m+n+2}) \in\mathcal{D}_{m+1|n+1} $, one has
\begin{gather}
\label{eq:deflrTidA}
 ({\bf l}\boxtimes^{A(m,n)}{\bf r})(T_{i,d})=
 \begin{cases}
  P_{\sigma_i \tr d}\circ \iota_d
   &
  \text{if $p_i\ne p_{i+1}$,} \\
  \iota_d\circ({\bf l}(T_{\tau_ {+.d}^{-1}(i)})\otimes{\rm{id}}_W)\circ P_d  
   &
  \text{if $p_i= p_{i+1}=0$,} \\
    \iota_d\circ({\rm{id}}_V \otimes{\bf r}(T_{\tau_ {-.d}^{-1}(i)}))\circ P_d  
   &
  \text{if $p_i= p_{i+1}=1$.}
 \end{cases}
\end{gather}

\end{proposition}

\begin{figure}
\setlength{\unitlength}{1mm}
\begin{picture}(30,70)(-40,0)

\put(-2,72){$p_i$}\put(8,72){$p_{i+1}$}\put(18,72){$p_{i+2}\,\,\in\{0,1\}$}
\put(-2,-3){$p_{i+2}$}\put(8,-3){$p_{i+1}$}\put(18,-3){$p_i$}

\put(-7,64){$d$}\put(-12,44){$i\sa d$}\put(-27,24){$(i+1)\sa i\sa d$}
\put(-32,4){$i\sa (i+1)\sa i\sa d$}

\put(-2,63.5){\line(0,1){3}}\put(22,63.5){\line(0,1){3}}
\put(-2,66.5){\line(1,0){24}}\put(-2,63.5){\line(1,0){24}}

\put(-2,43.5){\line(0,1){3}}\put(22,43.5){\line(0,1){3}}
\put(-2,46.5){\line(1,0){24}}\put(-2,43.5){\line(1,0){24}}

\put(-2,23.5){\line(0,1){3}}\put(22,23.5){\line(0,1){3}}
\put(-2,26.5){\line(1,0){24}}\put(-2,23.5){\line(1,0){24}}

\put(-2,3.5){\line(0,1){3}}\put(22,3.5){\line(0,1){3}}
\put(-2,6.5){\line(1,0){24}}\put(-2,3.5){\line(1,0){24}}

\put(0,60){\line(0,1){10}}\put(10,60){\line(0,1){10}}\put(20,60){\line(0,1){10}}
\put(0,50){\line(1,1){10}}\put(10,50){\line(-1,1){10}}\put(20,50){\line(0,1){10}}
\put(0,40){\line(0,1){10}}\put(10,40){\line(0,1){10}}\put(20,40){\line(0,1){10}}
\put(0,30){\line(0,1){10}}\put(10,30){\line(1,1){10}}\put(20,30){\line(-1,1){10}}
\put(0,20){\line(0,1){10}}\put(10,20){\line(0,1){10}}\put(20,20){\line(0,1){10}}
\put(0,10){\line(1,1){10}}\put(10,10){\line(-1,1){10}}\put(20,10){\line(0,1){10}}
\put(0,0){\line(0,1){10}}\put(10,0){\line(0,1){10}}\put(20,0){\line(0,1){10}}

\end{picture}
\setlength{\unitlength}{1mm}
\begin{picture}(30,70)(-80,0)

\put(-2,72){$p_i$}\put(8,72){$p_{i+1}$}\put(18,72){$p_{i+2}\,\,\in\{0,1\}$}
\put(-2,-3){$p_{i+2}$}\put(8,-3){$p_{i+1}$}\put(18,-3){$p_i$}

\put(-7,64){$d$}\put(-22,44){$(i+1)\sa d$}\put(-27,24){$i\sa (i+1)\sa d$}
\put(-42,4){$(i+1)\sa i\sa (i+1)\sa d$}

\put(-2,63.5){\line(0,1){3}}\put(22,63.5){\line(0,1){3}}
\put(-2,66.5){\line(1,0){24}}\put(-2,63.5){\line(1,0){24}}

\put(-2,43.5){\line(0,1){3}}\put(22,43.5){\line(0,1){3}}
\put(-2,46.5){\line(1,0){24}}\put(-2,43.5){\line(1,0){24}}

\put(-2,23.5){\line(0,1){3}}\put(22,23.5){\line(0,1){3}}
\put(-2,26.5){\line(1,0){24}}\put(-2,23.5){\line(1,0){24}}

\put(-2,3.5){\line(0,1){3}}\put(22,3.5){\line(0,1){3}}
\put(-2,6.5){\line(1,0){24}}\put(-2,3.5){\line(1,0){24}}

\put(0,60){\line(0,1){10}}\put(10,60){\line(0,1){10}}\put(20,60){\line(0,1){10}}
\put(0,50){\line(0,1){10}}\put(10,50){\line(1,1){10}}\put(20,50){\line(-1,1){10}}
\put(0,40){\line(0,1){10}}\put(10,40){\line(0,1){10}}\put(20,40){\line(0,1){10}}
\put(0,30){\line(1,1){10}}\put(10,30){\line(-1,1){10}}\put(20,30){\line(0,1){10}}
\put(0,20){\line(0,1){10}}\put(10,20){\line(0,1){10}}\put(20,20){\line(0,1){10}}
\put(0,10){\line(0,1){10}}\put(10,10){\line(1,1){10}}\put(20,10){\line(-1,1){10}}
\put(0,0){\line(0,1){10}}\put(10,0){\line(0,1){10}}\put(20,0){\line(0,1){10}}

\end{picture}
  \caption{Braid relation}
  \label{fig:braid2}
\end{figure}
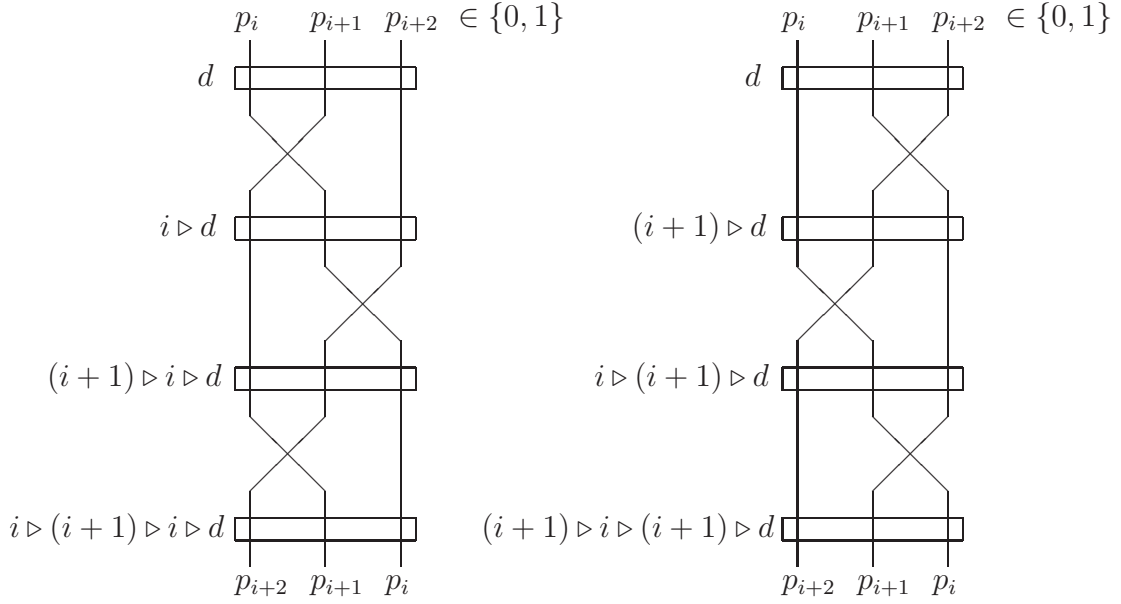

{\it {Proof.}} This can be checked directly. 
Refer to 
Figure~\ref{fig:braid2}. 
We explain by using an example. 
Denote $({\bf l}\boxtimes^{A(m,n)}{\bf r})(T_{i^\prime,d^\prime})$ 
by $S _{i^\prime,d^\prime}$ for any $d^\prime$ and $i^\prime$.
Let 
$d=(p_1,\ldots,p_{m+n+2})\in\mathcal{D}_{m+1|n+1} $
and $i\in\{1,\ldots,m+n\}$ and assume $p_i=p_{i+1}=0$
and $p_{i+2}=1$.
Let $d _1:=i\tr d(=\sigma_i \tr d)$, $d _2:=(i+1)\tr d _1$, $d _3:=i\tr d _2$,
$d _4:=(i+1)\tr d $, $d _5:=i\tr d _4$ and $d _6:=(i+1)\tr d _4$.
Then 
\begin{align}\label{eq:examofp}
d=d _1=&(p_1,\ldots, p_{i-1},0,0,1,p_{i+2}, \ldots,p_{m+n+2}), \\
d _2=d _4=&(p_1,\ldots, p_{i-1},0,1,0,p_{i+2}, \ldots,p_{m+n+2}), \\
d _3=d _5=d _6=&(p_1,\ldots, p_{i-1},1,0,0,p_{i+2}, \ldots,p_{m+n+2}). 
\end{align} Note that $\tau_ {+.d_5}=\sigma _i \sigma _{i+1} \tau_ {+.d}$.
Hence $\tau_ {+.d_5}^{-1}(i+1)=\tau_ {+.d}^{-1}(i)$.
Then we have 
$S _{i,d}=\iota _d\circ({\bf l}(T_{\tau_ {+.d}^{-1}(i)})\otimes{\rm{id}}_W)\circ P _d$,
$S _{i+1,d _1}=\iota_ {d _2}\circ P _d$,
$S _{i,d _2}=\iota_ {d _3}\circ P _{d _2}$,
$S _{i+1,d}=\iota_ {d _2}\circ P _d$,
$S _{i,d _4}=\iota_ {d _3}\circ P _{d _2}$
and
$S _{i+1,d _5}=\iota_ {d _3}\circ ({\bf l}(T_{\tau_ {+.d}^{-1}(i)})\otimes{\rm{id}}_W) \circ P _{d _3}$.
Hence we have $S _{i,d _2}S _{i+1,d _1}S _{i,d}=S _{i+1,d _5}S _{i,d _4}S _{i+1,d}
=\iota _{d _3}\circ({\bf l}(T_{\tau_ {+.d}^{-1}(i)})\otimes{\rm{id}}_W)\circ P _d$, as desired.
\hfill $\Box$
\newline\newline\par
For $\lambda\in \Lambda_{H_q(S_{m+1})}$
and $\mu\in \Lambda_{H_q(S_{n+1})}$, we denote 
$\rho_{H_q(S_{m+1}),\lambda}\boxtimes^{A(m,n)}\rho_{H_q(S_{n+1}),\mu}$
by $\rho^{A(m,n)}_{q;\lambda,\mu}$
and we denote $C_{V\otimes W}$, $P _d$, $\iota _d$ for $V=V_{H_q(S_{m+1}),\lambda}$ 
and $W=V_{H_q(S_{n+1}),\mu}$ by 
$C^{A(m,n)}_{q;\lambda,\mu}$, $P ^{\lambda,\mu} _d$, $\iota ^{\lambda,\mu} _d$
respectively.

\begin{theorem}\label{theorem:repThA}
Let $q\in\ndC$ and assume that 
\begin{align}\label{eq:assThA}
q P_{S_{m+1}}(q)P_{S_{n+1}}(q) \ne 0.
\end{align} Then the $\ndC$-algebra homomorphism
\begin{align}
&\bigoplus_{(\lambda,\mu)\in \Lambda_{H_q(S_{m+1})}
\times\Lambda_{H_q(S_{n+1})}} \rho^{A(m,n)}_{q;\lambda,\mu}
: & & \label{eq:bighomoftypeA} \\ & H_q(A(m,n))\rightarrow 
\bigoplus_{(\lambda,\mu)\in \Lambda_{H_q(S_{m+1})}
\times\Lambda_{H_q(S_{n+1})}} \End_\ndC (C^{A(m,n)}_{q;\lambda,\mu})
& & \nonumber 
\end{align} is an isomorphism. Further we have
\begin{equation}\label{eq:dimtypeA}
\dim H_q(A(m,n)) ={\frac {((m+n+2)!)^2} {(m+1)!(n+1)!}} .
\end{equation}
Moreover $H_q(A(m,n))$ is a
semisimple $\ndC$-algebra and a complete set of non-equivalent irreducible representations
of $H_q(A(m,n))$ is given by 
$\{\rho^{A(m,n)}_{q;\lambda,\mu} |(\lambda,\mu)\in \Lambda_{H_q(S_{m+1})}
\times\Lambda_{H_q(S_{n+1})}\}$. 

\end{theorem}

{\it {Proof.}} Define the $\ndC$-algebra homomorphism
$f_1:H_ q(S _{m+1})\otimes H_ q(S _{n+1})\rightarrow H_ q(A(m,n))$
by $f_1 (T_i\otimes 1)=T_{i,e}$ and $f_1 (1\otimes T_j)=T_{m+1+j,e}$.
Let $R _{\lambda,\mu}:=(\iota ^{\lambda,\mu} _{d_e}\circ P ^{\lambda,\mu} _{d_e})
\End_\ndC (C^{A(m,n)}_{q;\lambda,\mu})(\iota ^{\lambda,\mu} _{d_e}\circ 
P ^{\lambda,\mu} _{d_e})$.
Let $f_2$ denote the homomorphism of (\ref{eq:bighomoftypeA}).
It follows from (\ref{eq:assThA}) that $H_ q(S _{m+1})\otimes H_ q(S _{n+1})$ is
a semisimple $\ndC$-algebra. This implies 
$${\rm{Im}}(f_2\circ f_1)=\bigoplus_{(\lambda,\mu)\in \Lambda_{H_q(S_{m+1})}
\times\Lambda_{H_q(S_{n+1})}} R _{\lambda,\mu} .
$$ On the other hand, we have 
$$
\End_\ndC (C^{A(m,n)}_{q;\lambda,\mu})=\bigoplus _{d _1,d _2\in \mathcal{D}_{m+1|n+1}}
(\iota ^{\lambda,\mu} _{d _1}\circ P ^{\lambda,\mu} _{d_e})R _{\lambda,\mu}
(\iota ^{\lambda,\mu} _{d_e}\circ P ^{\lambda,\mu} _{d _2}) .
$$ Hence by (\ref{eq:deflrTidA}) we can easily see that $f_2$ is surjective.
In particular, we have
\begin{eqnarray*}
\lefteqn{\dim H_ q(A(m,n)) } \\
& \geq & \sum_{(\lambda,\mu)\in \Lambda_{H_q(S_{m+1})}
\times\Lambda_{H_q(S_{n+1})}}|\mathcal{D}_{m+1|n+1}| ^2\dim R _{\lambda,\mu} 
 \\
& = & |\mathcal{D}_{m+1|n+1}| ^2\sum_{(\lambda,\mu)\in \Lambda_{H_q(S_{m+1})}
\times\Lambda_{H_q(S_{n+1})}}\dim R _{\lambda,\mu} 
 \\ 
& =& 
({\frac {(m+n+2)!} {(m+1)!(n+1)!}})^2 \\
& & \sum_{(\lambda,\mu)\in \Lambda_{H_q(S_{m+1})} \times\Lambda_{H_q(S_{n+1})}}
 (\dim V_{H_q(S_{m+1}),\lambda})^2 (\dim V_{H_q(S_{n+1}),\mu})^2
 \\
& =& ({\frac {(m+n+2)!} {(m+1)!(n+1)!}})^2 \dim H_q(S_{m+1}) \dim H_q(S_{n+1}) 
\\  
& =& {\frac {((m+n+2)!)^2} {(m+1)!(n+1)!}}.
\end{eqnarray*}
Hence by \eqref{eq:spanHq3typeA},  
we have \eqref{eq:dimtypeA}.
Hence $f_2$ is an isomorphism.
Then the rest of the statement follows 
from well-known facts concerning semisimple algebras
(cf. \cite[(25.22) and (27.4)]{CR}). \hfill $\Box$

\subsection{Iwahori-Hecke type algebra
associated with the 
Lie superalgebra $\frak{osp}(2m+1|2n)$} \label{subsection:Bmn}

Let $m \in \ndN\cup\{0\}$ and $n \in \ndN$.
Ler $\ell :=m+n$.
For $1\leq i\leq \ell$, define 
${\hat{\sigma}}_i\in  S_\ell$ 
by 
${\hat{\sigma}}_i:=\sigma_i$
($1\leq i\leq \ell -1$).
and ${\hat{\sigma}}_\ell:={\rm{id}}$.

Let $W$ be the Coxeter groupoid 
associated with $(R,N,\is ,\sa )\in \srs $ such that
$N=\{1,2,\ldots,\ell\}$, $A=\mathcal{D}_{m|n}$, $i\sa d= {\hat{\sigma}}_i \tr d$,
$V_0=V^{(m+n)}_0$, $R_d^+=\{\varepsilon_i-\varepsilon_j ,\varepsilon_i-\varepsilon_j
| 1\leq i<j\leq \ell  \}\cup\{\varepsilon_i| 1\leq i\leq \ell  \}$,
$\al_{i,d}=\varepsilon_i-\varepsilon_{i+1}$ ($1\leq i\leq \ell -1$),
$\al_{\ell ,d}=\varepsilon_\ell$ and $\sigma_{i,d}={\widetilde{\sigma}}_{\al_{i,d}}$.
Denote $H_q(W)$ by $H_q(B(m,n))$.  
Then 
$H_q(B(m,n))$
is the $\ndC$-algebra (with $1$) generated by
\begin{align}
\{E_d\,|\,d\in\mathcal{D}_{m|n} \}
\cup\{T_{i,d}\,|\,1\leq i\leq \ell,\,d\in\mathcal{D}_{m|n} \}
	\label{eq:HBgen}
\end{align} 
and defined by the relations \eqref{eq:HIrel1-1}-\eqref{eq:HIrel4} and the relations
\begin{align}
	T_{\ell -1,{\hat{\sigma}}_{\ell -1}\tr d}T_{\ell , {\hat{\sigma}}_{\ell -1}\tr d}T_{\ell -1,d}T_{\ell , d}=&
	T_{\ell , d}T_{\ell -1,{\hat{\sigma}}_{\ell -1}\tr d}T_{\ell , {\hat{\sigma}}_{\ell -1}\tr d}T_{\ell -1,d}
    &  &\text{}
	\label{eq:HBrel9}\\
	T_{i,{\hat{\sigma}}_{i+1}{\hat{\sigma}}_i\tr d}
    T_{i+1,{\hat{\sigma}}_i\tr d}T_{i,d}=&
	T_{i+1,{\hat{\sigma}}_i{\hat{\sigma}}_{i+1}\tr d}
    T_{i,{\hat{\sigma}}_{i+1}\tr d}T_{i+1,d}&  
    &\text{if $1\leq i\leq \ell -1$,} 
    \label{eq:HBrel10}
	\\
	T_{i,{\hat{\sigma}}_j\tr d}T_{j,d}=&T_{j,{\hat{\sigma}}_i\tr d}T_{i,d}&
	&\text{if $|i-j|\geq 2$.}
	\label{eq:HBrel11}
\end{align}
\vspace{1cm}



\begin{figure}
  \begin{center}
\setlength{\unitlength}{1mm}
\begin{picture}(30,50)(-10,-8)

\put(-27,45){$d=(0,0,1)$}
\put(-27,25){$d=(0,1,0)$}
\put(-27,5){$d=(1,0,0)$}

\put(17,30){\line(0,1){10}}\put(14,34){$2$}
\put(2,10){\line(0,1){10}}\put(-1,14){$1$}

\put(2,45){\circle{3}}
\put(17,45){\circle{3}}\put(15.70,44.20){$\times$}
\put(32,45){\circle{3}}
\put(3.5, 45){\line(1,0){12}}
\put(18.5,45.5){\line(1,0){11}}\put(18.5,44.5){\line(1,0){11}}
\put(29,44){$\rangle$}

\put(2,25){\circle{3}}\put(0.70,24.20){$\times$}
\put(17,25){\circle{3}}\put(15.70,24.20){$\times$}
\put(32,25){\circle*{3}}
\put(3.5, 25){\line(1,0){12}}
\put(18.5,25.5){\line(1,0){11}}\put(18.5,24.5){\line(1,0){11}}
\put(29,24){$\rangle$}

\put(2,5){\circle{3}}\put(0.70,4.20){$\times$}
\put(17,5){\circle{3}}
\put(32,5){\circle*{3}}
\put(3.5, 5){\line(1,0){12}}
\put(18.5,5.5){\line(1,0){11}}\put(18.5,4.5){\line(1,0){11}}
\put(29,4){$\rangle$}

\end{picture}
  \end{center}
  \caption{Dynkin diagrams of the Lie superalgebra $B(1,2)$}
  \label{fig:B12}
\end{figure}
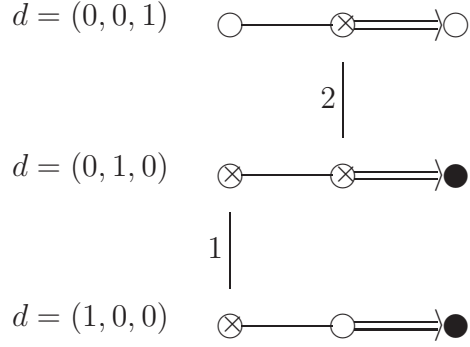

Recall $\rho$ and $d_e\in\mathcal{D}_{m|n}$ from 
Theorem~\ref{theorem:rep} and
\eqref{eq:dedo} respectively.
Then
$P_{d_e}\rho(e_{d_e}We_{d_e})\iota_{d_e}\subset\{(\sum_{i=1}^{m+n}z_i{\bf{E}}_{\sigma(i)i})|
\sigma\in S_{m+n}, z_i \in\{-1,1\}, \sigma(\{1,\ldots,m\})=\{1,\ldots,m\}\}$. Hence $|e_{d_e}We_{d_e}|\leq 2^{m+n}m!n!$
by Theorem~\ref{theorem:rep}, so 
$|W\setminus\{0\}|=|\mathcal{D}_{m|n}|^2|e_{d_e}W e_{d_e}|\leq {\frac {2^{m+n}((m+n)!)^2} {m!n!}}$.
Hence by \eqref{eq:spanHq3}, we conclude
\begin{align}\label{eq:spanHq3typeB}
\dim H_q(B(m,n))\leq {\frac {2^{m+n}((m+n)!)^2} {m!n!}}.
\end{align}

\begin{proposition}\label{proposition:Bmn}
Let $V _{\bf l}$ and $V _{\bf r}$ be finite dimensional $\ndC$-linear spaces, and
let ${\bf l}:H_q(W(B _m))\rightarrow
\End_\ndC (V _{\bf l})$ and 
${\bf r}:H_q(W(B _n))\rightarrow
\End_\ndC (V _{\bf r})$
be $\ndC$-algebra homomorphisms, i.e., representations.
Let ${\bf l}\otimes{\bf r}:H_q(W(B _m))\otimes 
H_q(W(B _n))\rightarrow 
\End_\ndC (V _{\bf l}\otimes V _{\bf r})$ denote the tensor representation of 
${\bf l}$ and ${\bf r}$ in the ordinary sense.
Let $C_{V _{\bf l}\otimes V _{\bf r};d}$ be copies of 
the $\ndC$-linear space
$V _{\bf l}\otimes V _{\bf r}$, indexed by $d\in\mathcal{D}_{m|n} $.
Let $C_{V _{\bf l}\otimes V _{\bf r}}:=
\oplus_{d\in\mathcal{D}_{m|n}}
C_{V _{\bf l}\otimes V _{\bf r};d}
$. 
Let $P_d:C_{V _{\bf l}\otimes V _{\bf r}}
\rightarrow C_{V\otimes W;d}$
and $\iota_d:C_{V\otimes W;d}\rightarrow 
C_{V _{\bf l}\otimes V _{\bf r}}$ denote 
the canonical
projection and the canonical 
inclusion map respectively.
Then there exists a unique $\ndC$-algebra homomorphism 
${\bf l}\boxtimes{\bf r}={\bf l}\boxtimes^{B(m,n)}{\bf r}:H_q(B(m,n))\rightarrow \End_\ndC 
(C_{V _{\bf l}\otimes V _{\bf r}})$ satisfying the following conditions:
\newline\par
{\rm{(i)}} For each $d\in\mathcal{D}_{m|n} $,
one has $({\bf l}\boxtimes{\bf r})(E_d)=\iota_d\circ P_d$, 
 \par
{\rm{(ii)}} For each $i\in\{1,\ldots,\ell =m+n \}$ and 
each $d=(p_1,\ldots,p_\ell) \in\mathcal{D}_{m|n} $, one has
\begin{gather}
\label{eq:deflrTidB}
 ({\bf l}\boxtimes{\bf r})(T_{i,d})=
 \begin{cases}
  P_{{\hat{\sigma}}_i \tr d}\circ \iota_d
   &
  \text{if $1\leq i\leq \ell -1$ and $p _i\ne p _{i+1}$,} \\
  \iota_d\circ({\bf l}(T_{\tau_ {+.d}^{-1}(i)})\otimes{\rm{id}}_{V _{\bf r}})\circ P_d  
   &
  \text{if $1\leq i\leq \ell -1$ and $p _i = p _{i+1} =0$,} \\
    \iota_d\circ({\rm{id}}_{V _{\bf l}} \otimes{\bf r}(T_{\tau_ {-.d}^{-1}(i)}))\circ P_d  
   &
  \text{if $1\leq i\leq \ell -1$ and $p _i = p _{i+1} =1$,} \\
    \iota_d\circ({\bf l}(T_ m )\otimes{\rm{id}}_{V _{\bf r}})\circ P_d
   &
  \text{if $i=\ell$ and $p _\ell =0$,} \\   
  \iota_d\circ({\rm{id}}_{V _{\bf l}}\otimes{\bf r}(T_ n ))\circ P_d
   &
  \text{if $i=\ell$ and $p _\ell =1$,} 
 \end{cases}
\end{gather} where $\tau_ {\pm.d}$ are the ones of \eqref{eq:taudplusminus}.

\end{proposition}
\vspace{1cm}
{\it {Proof.}} We can check out this directly
in a way similar to that for  Proof of  
Proposition~\ref{proposition:repA}.  \hfill $\Box$
\newline\newline\par
For $\lambda\in \Lambda_{H_q(W(B_m))}$
and $\mu\in \Lambda_{H_q(W(B_n))}$, we denote 
$\rho_{H_q(W(B_m)),\lambda}\boxtimes^{B(m,n)}\rho_{W(B_n)),\mu}$
by $\rho^{B(m,n)}_{q;\lambda,\mu}$
and we denote $C_{V\otimes W}$ for $V=V_{H_q(W(B_m)),\lambda}$ 
and $W=V_{H_q(W(B_n)),\mu}$ by 
$C^{B(m,n)}_{q;\lambda,\mu}$.

\begin{theorem}\label{theorem:repTypeB}
Let $q\in\ndC$ and assume that 
\begin{align}\label{eq:assThB}
q P_{W(B_m)}(q)P_{W(B_n)}(q) \ne 0.
\end{align} Then the $\ndC$-algebra homomorphism
\begin{align}
&\bigoplus_{(\lambda,\mu)\in \Lambda_{H_q(W(B_m))}
\times\Lambda_{H_q(W(B_n))}} \rho^{B(m,n)}_{q;\lambda,\mu}
: & & \label{eq:bighomoftypeB} \\ & H_q(B(m,n))\rightarrow 
\bigoplus_{(\lambda,\mu)\in \Lambda_{H_q(W(B_m))}
\times\Lambda_{H_q(W(B_n))}} \End_\ndC (C^{B(m,n)}_{q;\lambda,\mu})
& & \nonumber 
\end{align} is an isomorphism. Further we have
\begin{equation}\label{eq:dimtypeB}
\dim H_q(B(m,n)) ={\frac {2^{m+n}((m+n)!)^2} {m!n!}} .
\end{equation}
Moreover $H_q(B(m,n))$ is a
semisimple $\ndC$-algebra and a complete set of non-equivalent irreducible representations
of $H_q(B(m,n))$ is given by 
$\{\rho^{B(m,n)}_{q;\lambda,\mu} |(\lambda,\mu)\in \Lambda_{H_q(W(B_m))}
\times\Lambda_{H_q(W(B_n))}\}$. 
\end{theorem}
\vspace{1cm}
{\it {Proof.}} Let ${\bf l}:H_q(W(B _m))\rightarrow
\End_\ndC (V _{\bf l})$ and 
${\bf r}:H_q(W(B _n))\rightarrow
\End_\ndC (V _{\bf r})$ be irreducible representations.
Further, let  ${\bf l}\boxtimes{\bf r}:H_q(B(m,n))\rightarrow \End_\ndC
(C_{V _{\bf l}\otimes V _{\bf r}})$ be the representaion introduced in Proposition~\ref{proposition:Bmn}
for these ${\bf l}$ and ${\bf r}$.
By \eqref{eq:deflrTidB}, we can easily see that
\begin{align}\label{eq:onepfThB}
\forall d^\prime, \forall d^{\prime\prime}\in\mathcal{D}_{m|n}, \quad 
P_{d^\prime}\circ \iota_{d^{\prime\prime}} \in{\rm{Im}}({\bf l}\boxtimes{\bf r}).
\end{align}
Define the representation $f_1:H_q(W(B _m))\otimes H_q(W(B _n))
\rightarrow {\rm{End}}_\ndC(C_{V _{\bf l}\otimes V _{\bf r};d_e})$ by 
$f_1(T_i\otimes 1)= (P_{d_e}\circ\iota_{d_o})(({\bf l}\boxtimes{\bf r})(T_{n+i,d_o}))(P_{d_o}\circ\iota_{d_e})$
and $f_1(1\otimes T_j)= ({\bf l}\boxtimes{\bf r})(T_{m+j,d_e})$.
The condition \eqref{eq:assThB} implies that $f_1$ is an irreducible
representaion of $H_q(W(B _m))\otimes H_q(W(B _n))$.  
Moreover, using \eqref{eq:onepfThB}, we can easily see that ${\bf l}\boxtimes{\bf r}$ is an
irreducible representation of $H_q(B(m,n))$.
 
By the above argument, together with \eqref{eq:spanHq3typeB},  
in the same way as that for Proof of Theorem~\ref{theorem:repThA},
we can complete the proof of this theorem. \hfill $\Box$

\subsection{Iwahori-Hecke type algebra 
associated with the 
Lie superalgebra $\frak{osp}(2m|2n)$} \label{subsection:Dmn}

Let $m$, $n\in\ndN$. Define the set 
$\mathcal{D}^{CD}_{m|n}$ by
\begin{align}\label{eq:defDCD}
\mathcal{D}^{CD}_{m|n}:=&\{d^D|d=(p_1,\ldots,p_{m+n})\in\mathcal{D}_{m|n},\,p_{m+n}=0\} \\
& \cup\{d^C_+,d^C_-|d=(p_1,\ldots,p_{m+n})\in\mathcal{D}_{m|n},\,p_{m+n}=1\}, \nonumber
\end{align} so that
\begin{align}\label{eq:nmbDCD}
|\mathcal{D}^{CD}_{m|n}|={\frac {(m+n-1)!} {(m-1)!n!}}+2{\frac {(m+n-1)!} {m!(n-1)!}}
={\frac {(m+n-1)!(m+2n)} {m!n!}}.
\end{align}

Let $\ell:=m+n$ and $N=\{1,\ldots , \ell\}$.
Define the action $\tr$ of $F_2(N)$ on $\mathcal{D}^{CD}_{m|n}$ by
\begin{gather}
\label{eq:defCDtr}
i \tr a =
 \begin{cases}
  (\sigma_i \tr d)^D 
   &
  \text{if $a =d^D$, $1\leq i\leq \ell -2$ and $p_i\ne p_{i+1}$,} \\
  (\sigma_i \tr d)^C_+ 
   &
  \text{if $a =d^D$, $i= \ell -1$ and $p_i\ne p_{i+1}$,} \\
  (\sigma_{i-1} \tr d)^C_- 
   &
  \text{if $a =d^D$, $i= \ell $ and $p_{i-1}\ne p_i$,} \\
  (\sigma_i \tr d)^C_\pm 
   &
  \text{if $a =d^C_\pm$, $1\leq i\leq \ell -2$ and $p_i\ne p_{i+1}$,} \\
  (\sigma_i \tr d)^D 
   &
  \text{if $a =d^C_+$, $i= \ell -1$ and $p_i\ne p_{i+1}$,}  \\
  (\sigma_{i-1} \tr d)^D 
   &
  \text{if $a =d^C_-$, $i= \ell $ and $p_{i-1}\ne p_i$,}  \\
  a 
   &
  \text{otherwise.}  
 \end{cases}
\end{gather}
Now we define $R=(R,N,\mathcal{D}^{CD}_{m|n},\tr)\in\cR$ as follows. 
Let $N$ be as above.
Let $A=\mathcal{D}^{CD}_{m|n}$. Let $V_0=V^{(\ell)}_0$. 
Let $a=d^D$, $d^C_+$ or $d^C_-\in\mathcal{D}^{CD}_{m|n}$ 
with $d=(p_1,\ldots,p_{m+n})\in\mathcal{D}_{m|n}$.
Let $R_a^+$ be the subset of $V_0$ formed by the elements
$\varepsilon_i-\varepsilon_j$,
$\varepsilon_i+\varepsilon_j$, 
($1\leq i<j\leq \ell$) and $2\varepsilon_k$
($1\leq k\leq \ell$ and $p_k=1$). Define 
\begin{gather}
\label{eq:defCDpi}
\al_{i,a} :=
 \begin{cases}
 \varepsilon_i-\varepsilon_{i+1}
   &
  \text{if $a =d^D$ or $d^C_+$ and $1\leq i\leq \ell -1$,} \\
 \varepsilon_i-\varepsilon_{i+1}
   &
  \text{if $a =d^C_-$ and $1\leq i\leq \ell -2$,} \\
 \varepsilon_{\ell -1}+\varepsilon_\ell
   &
  \text{if $a =d^D$ and $i= \ell $,} \\
 2\varepsilon_\ell
   &
  \text{if $a =d^C_+$ and $i= \ell $,} \\
 -2\varepsilon_\ell 
   &
  \text{if $a =d^C_-$ and $i=\ell -1$,} \\
 \varepsilon_{\ell -1}+\varepsilon_\ell
   &
  \text{if $a =d^C_-$ and $i=\ell $.}
 \end{cases}
\end{gather} Define $\sigma_{i.a}:=
{\widetilde{\sigma}}_{\al_{i.a}}$.
Let $W$ be the Coxeter groupoid associated with
$R$. Recall $\rho$ and $d_e\in\mathcal{D}_{m|n}$ from 
Theorem~\ref{theorem:rep} and
\eqref{eq:dedo} respectively.
It is easy to show that 
\begin{align*} 
& \rho(e_{(d_e)^D}We_{(d_e)^D})=  & \\
&\{\sum_{j=1}^\ell z_j{\bf{E}}_{\sigma (j)j}|
\sigma\in S_\ell , z_j\in\{-1,1\},\prod_{j=n+1}^\ell z_j=1, 
\sigma(\{1,\ldots,n\})=\{1,\ldots,n\} \},  & \nonumber 
\end{align*} so $|e_{(d_e)^D}We_{(d_e)^D}|\leq m!n! 2^{\ell -1}$
by
Theorem~\ref{theorem:rep}.
Hence $|W\setminus\{0\}|\leq |\mathcal{D}^{CD}_{m|n}|^2m!n! 2^{\ell -1}$. 
Denote $H_q(W)$ by $H_q(\frak{osp}(2m|2n))$. 
By \eqref{eq:spanHq3} and \eqref{eq:nmbDCD}, we have
\begin{align}\label{eq:spanHq3typeCD}
\dim H_q(\frak{osp}(2m|2n))\leq {\frac {2^{m+n-1}((m+n-1)!(m+2n))^2} {m!n!}}.
\end{align}

Recall that 
$H_q(\frak{osp}(2m|2n))$
is the $\ndC$-algebra (with $1$) generated by
\begin{align}
\{E_a\,|\,a\in\mathcal{D}^{CD}_{m|n} \}
\cup\{T_{i,a}\,|\,1\leq i\leq m+n,\,a\in\mathcal{D}^{CD}_{m|n} \}
	\label{eq:HCDgen}
\end{align}
and defined by the relations \eqref{eq:HIrel1-1}-\eqref{eq:HIrel4}
and the relations 
\begin{align}
	&(T_{i,a}T_{j,a})^2= 
	(T_{j,a}T_{i,a})^2
    \,\,\makebox[0pt][l]{if $a =d^C_\pm$, $p_{\ell -1}= p_\ell $ and $i=\ell -1$, $j=\ell $,} 
    &  &
	\label{eq:HCDrel9}
    \\
	&T_{i,a}T_{j,a}= 
	T_{j,a}T_{i,a}  
\,\,\makebox[0pt][l]{if $a =d^D$, $p_{\ell -1}= p_\ell $ and $i=\ell -1$, $j=\ell $,} 
&  
	& 
	\label{eq:HCDrel10}
    \\
	&T_{i,ji\tr d}T_{j,i\tr d}T_{i,d}=
	T_{j,ij\tr d}T_{i,j\tr d}T_{j,d} 
\,\,\makebox[0pt][l]{if $p_{\ell -1}\ne p_\ell $ and $i=\ell -1$, $j=\ell $,}&  &  
    \label{eq:HCDrel11}
    \\
	&T_{i,ji\tr d}T_{j,i\tr d}T_{i,d}=
	T_{j,ij\tr d}T_{i,j\tr d}T_{j,d}&  &  
    \text{if $1\leq i\leq \ell -3$, $j=i+1$,} 
    \label{eq:HCDrel12}
    \\
	&T_{i,ji\tr d}T_{j,i\tr d}T_{i,d}=
	T_{j,ij\tr d}T_{i,j\tr d}T_{j,d}&  &  
    \text{if $a =d^C_+$ and $i=\ell -2$, $j=\ell -1$}, 
    \label{eq:HCDrel13}
    \\ 
    &T_{i,ji\tr d}T_{j,i\tr d}T_{i,d}=
	T_{j,ij\tr d}T_{i,j\tr d}T_{j,d}&  &  
    \text{if $a =d^C_-$ and $i=\ell -2$, $j=\ell $,}
    \label{eq:HCDrel14} 
    \\
	&T_{i,ji\tr d}T_{j,i\tr d}T_{i,d}=
	T_{j,ij\tr d}T_{i,j\tr d}T_{j,d}
\,\,\makebox[0pt][l]{if $a =d^D$ and $i=\ell -2$, $j\in\{\ell -1, \ell \}$,}&  &  
    \label{eq:HCDrel15}
    \\
	&T_{j,i\tr d}T_{i,d}=
	T_{i,j\tr d}T_{j,d}    
	\,\,\makebox[0pt][l]{if $i<j$ holds and $i$, $j$ 
    are not the ones in \eqref{eq:HCDrel14} - \eqref{eq:HCDrel15}.}
	& & 
	\label{eq:HCDrel16}
\end{align}

\begin{figure}
  \begin{center}
\setlength{\unitlength}{1mm}
\begin{picture}(80,50)(10,-8)

\put(-27,45){$(1,0,0,0)^D$}
\put(-27,25){$(0,1,0,0)^D$}\put(53,30){$(0,0,0,1)^C_-$}
\put(-27,5){$(0,0,1,0)^D$}\put(53,5){$(0,0,0,1)^C_+$}

\put(2,30){\line(0,1){10}}\put(-1,34){$1$}
\put(17,10){\line(0,1){10}}\put(14,14){$2$}

\put(39.5,10){\line(1,1){12}}\put(45,19){$4$}
\put(39.5,1){\line(1,0){12}}\put(45,2){$3$}

\put(1, 48){$1$}\put(16, 48){$2$}\put(35,40){$3$}\put(35,48){$4$}
\put(2,45){\circle{3}}\put(0.70,44.20){$\times$}
\put(17,45){\circle{3}}
\put(32,41){\circle{3}}\put(32,49){\circle{3}}
\put(3.5, 45){\line(1,0){12}}
\put(18.5,44){\line(2,-1){6}}\put(24.5,41){\line(1,0){6}}
\put(18.5,46){\line(2,1){6}}\put(24.5,49){\line(1,0){6}}

\put(1,19){$1$}\put(16,28){$2$}\put(35,20){$3$}\put(35,28){$4$}
\put(2,25){\circle{3}}\put(0.70,24.20){$\times$}
\put(17,25){\circle{3}}\put(15.70,24.20){$\times$}
\put(32,21){\circle{3}}\put(32,29){\circle{3}}
\put(3.5,25){\line(1,0){12}}
\put(18.5,24){\line(2,-1){6}}\put(24.5,21){\line(1,0){6}}
\put(18.5,26){\line(2,1){6}}\put(24.5,29){\line(1,0){6}}

\put(1,-1){$1$}\put(16,-1){$2$}\put(35,0){$3$}\put(35,8){$4$}
\put(2,5){\circle{3}}
\put(17,5){\circle{3}}\put(15.70,4.20){$\times$}
\put(32,1){\circle{3}}\put(32,9){\circle{3}}
\put(30.70,0.20){$\times$}\put(30.70,8.20){$\times$}
\put(3.5,5){\line(1,0){12}}
\put(18.5,4){\line(2,-1){6}}\put(24.5,1){\line(1,0){6}}
\put(18.5,6){\line(2,1){6}}\put(24.5,9){\line(1,0){6}}
\put(31.5,2.5){\line(0,1){5}}\put(32.5,2.5){\line(0,1){5}}

\put(61, 19){$1$}\put(76,19){$2$}\put(91,19){$4$}\put(106,19){$3$}
\put(62,25){\circle{3}}
\put(77,25){\circle{3}}
\put(92,25){\circle{3}}\put(90.70,24.20){$\times$}
\put(107,25){\circle{3}}
\put(63.5, 25){\line(1,0){12}}\put(78.5,25){\line(1,0){12}}
\put(93.5,24.5){\line(1,0){12}}\put(93.5,25.5){\line(1,0){12}}

\put(61, -5){$1$}\put(76,-5){$2$}\put(91,-5){$3$}\put(106,-5){$4$}
\put(62,1){\circle{3}}
\put(77,1){\circle{3}}
\put(92,1){\circle{3}}\put(90.70,0.20){$\times$}
\put(107,1){\circle{3}}
\put(63.5, 1){\line(1,0){12}}\put(78.5,1){\line(1,0){12}}
\put(93.5,0.5){\line(1,0){12}}\put(93.5,1.5){\line(1,0){12}}

\end{picture}
  \end{center}
  \caption{Dynkin diagrams of the Lie superalgebra $D(3,1)$}
  \label{fig:D31}
\end{figure}
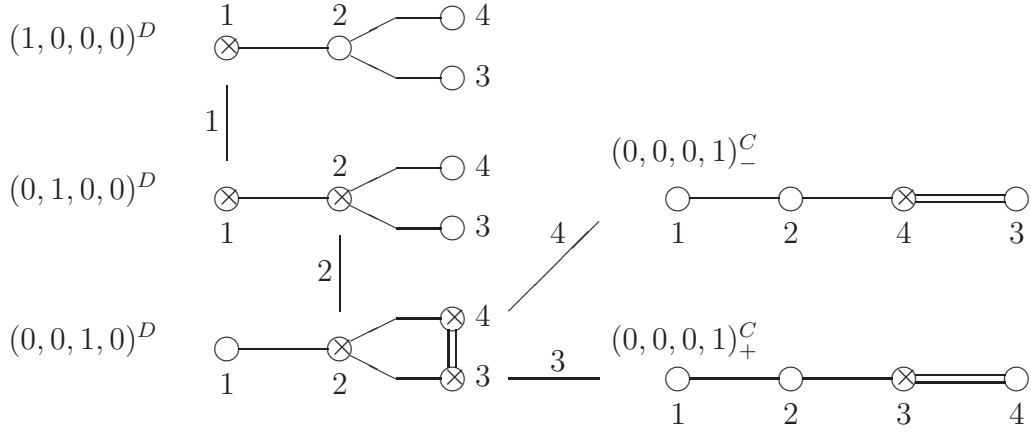

Recall that $W(C_k)=W(B_k)$ and 
$H_q(W(C_k))=H_q(W(B_k))$.

\begin{proposition}
Let $V _{\bf l}$ and $V _{\bf r}$ be finite dimensional $\ndC$-linear spaces, and
let ${\bf l}:H_q(W(D _m))\rightarrow
\End_\ndC (V _{\bf l})$ and 
${\bf r}:H_q(W(C _n))\rightarrow
\End_\ndC (V _{\bf r})$
be $\ndC$-algebra homomorphisms, i.e., representations.
Let ${\bf l}\otimes{\bf r}:H_q(W(D _m))\otimes 
H_q(W(C _n))\rightarrow 
\End_\ndC (V _{\bf l}\otimes V _{\bf r})$ denote the tensor representation of 
${\bf l}$ and ${\bf r}$ in the ordinary sense.
Let $C_{V _{\bf l}\otimes V _{\bf r};a}$ be copies of 
the $\ndC$-linear space
$V _{\bf l}\otimes V _{\bf r}$, indexed by $a\in\mathcal{D}^{CD}_{m|n} $.
Let $C_{V _{\bf l}\otimes V _{\bf r}}:=
\oplus_{d\in\mathcal{D}^{CD}_{m|n}}
C_{V _{\bf l}\otimes V _{\bf r};d}
$. 
Let $P_a:C_{V _{\bf l}\otimes V _{\bf r}}
\rightarrow C_{V\otimes W;a}$
and $\iota_a:C_{V\otimes W;a}\rightarrow 
C_{V _{\bf l}\otimes V _{\bf r}}$ denote 
the canonical
projection and the canonical 
inclusion map respectively.
Then there exists a unique $\ndC$-algebra homomorphism 
${\bf l}\boxtimes{\bf r}={\bf l}\boxtimes^{CD}{\bf r}:H_q(\frak{osp}(2m|2n))\rightarrow \End_\ndC 
(C_{V _{\bf l}\otimes V _{\bf r}})$ satisfying the following conditions:
\newline\par
{\rm{(i)}} For each $a\in\mathcal{D}^{CD}_{m|n} $,
one has $({\bf l}\boxtimes{\bf r})(E_a)=\iota_a\circ P_a$. 
 \par
{\rm{(ii)}} For each $i\in\{1,\ldots,\ell =m+n \}$ and 
each $a\in\mathcal{D}^{CD}_{m|n} $ with
$d=(p_1,\ldots,p_\ell) \in\mathcal{D}_{m|n} $ such that
$a=d^C_+$, $d^C_-$ or $d^D$, one has
\begin{gather}
\label{eq:deflrTiaCD}
 ({\bf l}\boxtimes{\bf r})(T_{i,a})=
 \begin{cases}
  P_{i \tr a}\circ \iota_a
   &
  \text{if $1\leq i\leq \ell  $ and $i \tr a\ne a$,} \\
  \iota_a
\circ({\bf l}(T_{\tau_ {+.d}^{-1}(i)})\otimes{\rm{id}}_{V _{\bf r}})\circ P_a  
   &
  \text{if $1\leq i\leq \ell  -1$ and $p _i = p _{i+1} =0$,} \\
    \iota_a\circ({\bf l}(T_ m )\otimes{\rm{id}}_{V _{\bf r}})\circ P_a
   &
  \text{if $i=\ell $ and $p _\ell  =0$,} \\
    \iota_a\circ({\rm{id}}_{V _{\bf l}} \otimes{\bf r}(T_{\tau_ {-.d}^{-1}(i)}))\circ P_a  
   &
  \text{if $1\leq i\leq \ell  -2$ and $p _i = p _{i+1} =1$,} \\   
  \iota_a\circ({\rm{id}}_{V _{\bf l}}\otimes{\bf r}(T_ {n-1} ))\circ P_a
   &
  \text{if $i=\ell -1$ and $p _ {\ell  -1}=p _\ell =1$, $a=d^C_+$.} \\   
  \iota_a\circ({\rm{id}}_{V _{\bf l}}\otimes{\bf r}(T_ n ))\circ P_a
   &
  \text{if $i=\ell  -1$ and $p _\ell =1$, $a=d^C_-$.} \\   
  \iota_a\circ({\rm{id}}_{V _{\bf l}}\otimes{\bf r}(T_ n ))\circ P_a
   &
  \text{if $i=\ell $ and $p _\ell  =1$, $a=d^C_+$.} \\   
  \iota_a\circ({\rm{id}}_{V _{\bf l}}\otimes{\bf r}(T_ {n-1} ))\circ P_a
   &
  \text{if $i=\ell $ and $p _ {\ell  -1}=p _\ell  =1$, $a=d^C_-$,} 
 \end{cases}
\end{gather} where $\tau_ {\pm.d}\in S_{m+n}$ are the ones of 
\eqref{eq:taudplusminus}.

\end{proposition}

{\it{Proof.}} We can check out this directly
in a way similar to that for  Proof of  
Proposition~\ref{proposition:repA}. \hfill $\Box$
\newline\newline\par
For $\lambda\in \Lambda_{H_q(W(D_m))}$
and $\mu\in \Lambda_{H_q(W(C_n))}$, we denote 
$\rho_{H_q(W(D_m)),\lambda}\boxtimes^{CD}\rho_{W(C_n)),\mu}$
by $\rho^{CD}_{q;\lambda,\mu}$
and we denote $C_{V\otimes W}$ for $V=V_{H_q(W(D_m)),\lambda}$ 
and $W=V_{H_q(W(C_n)),\mu}$ by 
$C^{CD}_{q;\lambda,\mu}$.

\begin{theorem}\label{theorem:repTypeCD}
Let $q\in\ndC$ and assume that 
\begin{align}\label{eq:assThCD}
q P_{W(D_m)}(q)P_{W(C_n)}(q) \ne 0.
\end{align} Then the $\ndC$-algebra homomorphism
\begin{align}
&\bigoplus_{(\lambda,\mu)\in \Lambda_{H_q(W(D_m))}
\times\Lambda_{H_q(W(C_n))}} \rho^{CD}_{q;\lambda,\mu}
: & & \label{eq:bighomoftypeCD} \\ & H_q(\frak{osp}(2m|2n))\rightarrow 
\bigoplus_{(\lambda,\mu)\in \Lambda_{H_q(W(D_m))}
\times\Lambda_{H_q(W(C_n))}} \End_\ndC (C^{CD}_{q;\lambda,\mu})
& & \nonumber 
\end{align} is an isomorphism. Further we have
\begin{equation}\label{eq:dimtypeCD}
\dim H_q(\frak{osp}(2m|2n)) ={\frac {2^{m+n-1}((m+n-1)!(m+2n))^2} {m!n!}}.
\end{equation}
Moreover $H_q(\frak{osp}(2m|2n))$ is a
semisimple $\ndC$-algebra and a complete set of non-equivalent irreducible representations
of $H_q(\frak{osp}(2m|2n))$ is given by 
$\{\rho^{CD}_{q;\lambda,\mu} |(\lambda,\mu)\in \Lambda_{H_q(W(D_m))}
\times\Lambda_{H_q(W(C_n))}\}$. 
\end{theorem}

{\it{Proof.}} Note that $W(D_m)\times W(C_n)=2^{m+n-1}m!n!$. Then we can prove
this theorem in the same way as that for Proof of Theorem~\ref{theorem:repTypeB}.
\hfill $\Box$

\begin{remark}
Now, by \eqref{eq:evoodptofgl}, \eqref{eq:evenpartofosp}
and Theorems~\ref{theorem:repThA}, \ref{theorem:repTypeB} and 
\ref{theorem:repTypeCD}, 
it has turned out that if $q$ is non-zero and not any primitive 
root of unity, then as a $\ndC$-algebra,
$H_q(\mathfrak{g})=H_q(W)$ introduced in this section for 
the Lie superalgebra $\mathfrak{g}=A(m,n)$ or $\mathfrak{osp}(m|2n)$
is very similar to the Iwahori-Hecke algebra $H_q(W_0)$ 
associated with the Weyl group $W_0$ of the Lie algebra $\mathfrak{g}(0)$
given as the even part of $\mathfrak{g}$.   
\end{remark}

\begin{remark} Assume $q$ to be an element of $\ndC$ transcendental over $\ndQ$. 
Then the $\ndZ$-subalgebra (with identity) of $\ndC$ generated by $q$
can also be regarded as the polynomial ring $\ndZ [q]$ in the variable $q$ over $\ndZ$.
Let $W$ be one of the Coxeter groupoids treated in Subsections~\ref{subsection:Amn},
\ref{subsection:Bmn} and \ref{subsection:Dmn}. By Lemma~\ref{lemma:spanHq}
and \eqref{eq:dimtypeA}, \eqref{eq:dimtypeB}, \eqref{eq:dimtypeCD}, one see that
$\{f(w)|w\in W\setminus\{0\}\}$ is a $\ndC$-basis of $H_q(W)$, that is, 
$H_q(W)=\oplus_{w\in W\setminus\{0\}}\ndC f(w)$.
Define $H_{\ndZ [q],q}(W)$ to be the $\ndZ [q]$-submodule 
generated by $\{f(w)|w\in W\setminus\{0\}\}$.
Using Theorem~\ref{theorem:matsu} and Corollary~\ref{corollary:matsu},
one see that $H_{\ndZ [q],q}(W)$ is also a $\ndZ [q]$-subalgebra of $H_q(W)$.
Let $\ndA$ be any commutative ring (with identity). Let 
$\zeta$ be any element of $\ndA$. Regard $\ndA$ as a $\ndZ [q]$-algebra via the
$\ndZ$-algebra homomorphism $\ndZ [q]\to \ndA$ that takes $q$ to $\zeta$.
Let $H_{\ndA ,\zeta}(W)$ be the $\ndA$-algebra (with identity) 
defined by 
$H_{\ndA ,\zeta}(W):=H_{\ndZ [q],q}(W)\otimes _{\ndZ [q]} \ndA$. For $X\in H_{\ndZ [q],q}(W)$,
we also denote the element $X\otimes 1$ of $H_{\ndA ,\zeta}(W)$ by $X$.
Then
$H_{\ndA ,\zeta}(W)$ is a free $\ndA$-module with an $\ndA$-basis 
$\{f(w)|w\in W\setminus\{0\}\}$, that is,
\begin{equation}
{\rm{rank}} _\ndA H_{\ndA ,\zeta}(W) = |W|-1.
\end{equation}
Using Theorem~\ref{theorem:matsu} and Corollary~\ref{corollary:matsu} again,
one see that $H_{\ndA ,\zeta}(W)$ can also be defined by the same generators
as \eqref{eq:genHalgebroid} and the same relations
as \eqref{eq:HIrel1-1}-\eqref{eq:HIrel6} with $\zeta$ in place of $q$.

The same properties as above seem to be true for many Coxeter groupoids,
which might be able to be proved in a way similar to that of the proof of \cite[Proposition~3.3]{Lu}.
\end{remark}

\providecommand{\bysame}{\leavevmode\hbox to3em{\hrulefill}\thinspace}
\providecommand{\MR}{\relax\ifhmode\unskip\space\fi MR }
\providecommand{\MRhref}[2]{%
  \href{http://www.ams.org/mathscinet-getitem?mr=#1}{#2}
}
\providecommand{\href}[2]{#2}

{\small \sc Hiroyuki Yamane,
Department of Pure and Applied Mathematics,
Graduate School of Information Science
and Technology, Osaka University, Toyonaka 560-0043,
Japan}

{\small \textit{E-mail address:} \texttt{yamane@ist.osaka-u.ac.jp}}

\end{document}